\documentclass[twoside,11pt]{article}

\usepackage{amsfonts}

\newtheorem{theorem}{Theorem}[section]
\newtheorem{lemma}[theorem]{Lemma}
\newtheorem{statement}[theorem]{Statement}
\newtheorem{corollary}[theorem]{Corollary}

\long\gdef\boxit#1{\begingroup\vbox{\hrule\hbox{\vrule\kern3pt
      \vbox{\kern3pt#1\kern3pt}\kern3pt\vrule}\hrule}\endgroup}

\def\qed{ \ \vrule width.2cm height.2cm depth0cm\smallskip}
\def\qedt{ \ \vrule width.2cm height.2cm depth0cm}

\def\Rset{{\mathbb R}}
\def\Zset{{\mathbb Z}}
\def\Nset{{\mathbb N}}
\def\Bscr{{\cal B}}
\def\Cscr{{\cal C}}
\def\Dscr{{\cal D}}

\def\Vscr{{\cal V}}

\def\div{\mbox{div}}
\def\max{{\rm max}}
\def\min{{\rm min}}
\def\supp{\mbox{supp}}
\def\lxipi{\ell^\xi_\pi}

\def\rdist{\mbox{r-dist}}
\def\tilde{\widetilde}
\def\bar{\overline}

\def\spitem{\vspace{-5pt}\item}

\newenvironment{proof}{\noindent {\bf Proof.\/}}{$\qed$\vskip 0.1in}
\newcommand{\Xcomment}[1]{}
\newcommand{\Remark}{{\bf Remark.\ }}

\newcommand{\refeq}[1]{(\ref{eq:#1})}

\newenvironment{myitem}{\refstepcounter{equation}\begin{enumerate}%
\item[(\arabic{equation})]$\quad$}{\end{enumerate}}

\makeatletter
\renewcommand{\section}{\@startsection{section}{1}{0pt}%
{-3.5ex plus -1ex minus -.2ex}{2.3ex plus .2ex}%
{\normalfont\Large}}

\renewcommand{\subsection}{\@startsection{subsection}{2}{0pt}%
{-3.25ex plus -1ex minus -.2ex}{1.5ex plus .2ex}%
{\normalfont\large\bf}}

\makeatother

\begin{document}
\begin{titlepage}

\def\thepage {} % Kill pagenumbering

\begin{flushleft}
\large Andrew V. Goldberg\footnote{A.V. Goldberg: InterTrust
Technologies Corp., 4750 Patrick Henry Drive, Santa Clara, CA 95054
email: avg@acm.org. Part of this research was done while this author
was at NEC Research Institute, Inc., Princeton, NJ.}
and Alexander V. Karzanov\footnote{A.V. Karzanov: Corresponding author.
Institute for System Analysis of the Russian Academy of Sciences, 9,
Prospect 60 Let Oktyabrya, 117312 Moscow, Russia, email:
sasha@cs.isa.ac.ru. This author was supported in part by a grant from
the Russian Foundation of Basic Research.}
  \end{flushleft}

\bigskip
\begin{flushleft}
\bf\Large Maximum skew-symmetric flows and matchings
 \end{flushleft}

\medskip
\begin{flushleft}
December 2003
\end{flushleft}
%Received: June 1999 / Accepted: December 2003

\bigskip
\noindent
{\bf Abstract.}
The maximum integer skew-symmetric flow problem (MSFP)
generalizes both the maximum flow and maximum matching
problems. It was introduced by Tutte~\cite{tut-67} in terms of
self-conjugate flows in antisymmetrical digraphs. He showed that for
these objects there are natural analogs of classical theoretical
results on usual network flows, such as the flow decomposition,
augmenting path, and max-flow min-cut theorems. We give unified and
shorter proofs for those theoretical results.

We then extend to MSFP the shortest
augmenting path method of Edmonds and Karp~\cite{EK-72} and the
blocking flow method of Dinits~\cite{din-70}, obtaining algorithms
with similar time bounds in general case. Moreover, in the cases of
unit arc capacities and unit ``node capacities'' the blocking
skew-symmetric flow algorithm has time bounds similar to those
established in~\cite{ET-75,kar-73-2} for Dinits' algorithm.
In particular, this implies an algorithm for finding a
maximum matching in a nonbipartite graph in $O(\sqrt{n}m)$ time,
which matches the time bound for the algorithm of Micali and
Vazirani~\cite{MV-80}.
Finally, extending a clique compression technique of Feder and
Motwani~\cite{FM-91} to particular skew-symmetric graphs, we
speed up the implied maximum matching algorithm to run
in $O(\sqrt{n}m\log(n^2/m)/\log{n})$ time, improving the
best known bound for dense nonbipartite graphs.

Also other theoretical and algorithmic results on skew-symmetric
flows and their applications are presented.



\end{titlepage}

\baselineskip 15pt

\section{\Large Introduction}

By a {\em skew-symmetric flow} we mean a flow in a
skew-symmetric directed graph which takes equal values on any pair of
``skew-symmetric'' arcs. This is a synonym of Tutte's
{\em self-conjugate flow in an antisymmetrical digraph}~\cite{tut-67}.
This paper is devoted to the maximum integer skew-symmetric flow
problem, or, briefly, the {\em maximum IS-flow problem}. We study
combinatorial properties of this problem and develop fast algorithms
for it.

A well-known fact~\cite{FF-62} is that the bipartite matching problem
can be viewed as a special case of the maximum flow problem.
The combinatorial structure of nonbipartite matchings revealed by
Edmonds~\cite{edm-65} involves blossoms and is
more complicated than the structure of flows.
This phenomenon explains, to some extent, why
general matching algorithms are typically more intricate relative
to flow algorithms. The maximum IS-flow problem is a generalization
of both the maximum flow and maximum matching (or b-matching) problems.
Moreover,
this generalization appears to be well-grounded for two reasons.
First, the basic combinatorial and linear programming theorems for
usual flows have natural counterparts for IS-flows. Second, when
solving problems on IS-flows, one can use intuition, ideas and
technical tools well-understood for usual flows, so that the implied
algorithms for matchings become more comprehensible.

As the maximum flow problem is related to certain path
problems, the maximum IS-flow problem is related to certain problems on
so-called {\em regular paths} in skew-symmetric graphs.
We use some theoretical and algorithmic results on the {\em regular
reachability} and {\em shortest regular path} problems
from~\cite{GK-96}.

Tutte~\cite{tut-67} originated a mini-theory of IS-flows (in our terms)
to bridge theoretical results on matchings and their
generalizations (b-factors, b-matchings, degree-constrained subgraphs,
Gallai's track packings, and etc.) and results on usual flows. This
theory parallels Ford and Fulkerson's flow theory~\cite{FF-62} and
includes as basic results the decomposition, augmenting path, and
max-flow min-cut theorems. Subsequently, some of those results were
re-examined in different, but equivalent, terms by other authors, e.g.,
in~\cite{blu-90,GK-95,KS-93}.

Recall that the flow decomposition theorem says that a flow can be
decomposed into a collection of source-to-sink paths and cycles. The
augmenting path theorem says that a flow is maximum if and only if it
admits no augmenting path. The max-flow min-cut theorem says that the
maximum flow value is equal to the minimum cut capacity. Their
skew-symmetric analogs are, respectively, that an IS-flow can be
decomposed into a collection of pairs of symmetric source-to-sink paths
and pairs of symmetric cycles, that an IS-flow is maximum if and only
if it admits no regular augmenting path, and that the maximum IS-flow
value is equal to the minimum {\em odd-barrier} capacity. We give
unified and shorter proofs for these skew-symmetric flow theorems.

There is a relationship between skew-symmetric flows and {\em
bidirected flows} introduced by Edmonds and Johnson~\cite{EJ-70}
in their combinatorial study of a natural class of integer linear
programs generalizing usual flow and matching problems. In
particular, they established a linear programming description for
integer bidirected flows. We finish the theoretical part by showing how
to obtain a linear programming description for maximum IS-flows
directly, using the max-IS-flow min-barrier theorem.

The second, larger, part of this paper is devoted to efficient
methods to solve the maximum IS-flow problem (briefly, {\em MSFP}) in
general and special cases, based on the theoretical ground given in the
first part. First of all we explain how to adapt the idea of
Anstee's elegant methods~\cite{ans-85,ans-87} for b-matchings in which
standard flow algorithms are used to construct an optimal half-integer
solution and then, after rounding, the ``good pre-solution''
is transformed into an optimal b-matching by solving $O(n)$
certain path problems.
We devise an $O(M(n,m)+nm)$-time algorithm for MSFP in a similar
fashion, using a regular reachability algorithm with linear
complexity to improve a good pre-solution.
Hereinafter $n$ and $m$ denote the numbers of nodes and arcs of the
input graph, respectively, and $M(n,m)$ is the time needed to find an
integer maximum flow in a usual network with $n$ nodes and $m$ arcs.
Without loss of generality, we assume $n=O(m)$.

The next approach is the core of this paper. The purpose is to
extend to MSFP the well-known shortest augmenting path algorithm of
Edmonds and Karp~\cite{EK-72} with complexity $O(nm^2)$, and its
improved version, the blocking flow algorithm of
Dinits~\cite{din-70} with complexity $O(n^2m)$, so as to preserve the
complexity bounds. Recall that the blocking flow algorithm consists
of $O(n)$ {\em phases}, each finding a blocking flow in a
layered network representing the union of currently shortest
augmenting paths. We introduce concepts of
shortest blocking and totally blocking IS-flows and show that an
optimal solution to MSFP is also constructed in $O(n)$ phases,
each finding a shortest blocking IS-flow in the residual skew-symmetric
network. In its turn a phase is reduced to finding a totally blocking
IS-flow in an {\em acyclic} (though not necessarily layered)
skew-symmetric network.

The crucial point is to perform the latter task in time comparable with
the phase time in Dinits' algorithm (which is $O(nm)$ in general case).
We reduce it to a certain auxiliary problem in a usual acyclic
digraph. A fast algorithm for this problem provides the desired time
bound for a phase.

The complexity of our blocking IS-flow algorithm
remains comparable with that of Dinits' algorithm in important
special cases where both the number of phases and the phase time
significantly decrease. More precisely, Dinits' algorithm applied
to the maximum matching problem in a bipartite graph runs in
$O(\sqrt{n}m)$ time~\cite{HK-73,kar-73}. Extending that result,
it was shown in~\cite{ET-75,kar-73-2} that for arbitrary
nonnegative integer capacities, Dinits' algorithm has
$O(\min\{n,\sqrt{\Delta}\})$ phases and each phase runs in
$O(\min\{nm,m+\Delta\})$ time, where $\Delta$ is the sum of transit
capacities of inner nodes. Here the transit capacity of a node
(briefly, the {\em node capacity}) is the maximum flow value that can
be pushed through this node. We show that both bounds remain valid for
the blocking IS-flow algorithm.

When the network has unit arc capacities (resp. unit inner node
capacities), the number of phases turns into $O(\sqrt{m})$ (resp.
$O(\sqrt{n})$); in both cases the phase time turns into $O(m)$.
The crucial auxiliary problem (that we are able to solve in linear
time for unit arc capacities) becomes the following {\em maximal
balanced path-set problem}:
  \begin{itemize}
\item[{\em MBP}:] {\sl Given an acyclic digraph in which one sink and
an even set of sources partitioned into pairs are distinguished,
find an (inclusion-wise) maximal set of pairwise arc-disjoint paths
from sources to the sink such that for each pair $\{z,z'\}$ of sources
in the partition, the number of paths from $z$ is equal to that from
$z'$.}
  \end{itemize}

As a consequence, the implied algorithm solves the maximum matching
problem in a general graph in the same time, $O(\sqrt{n}m)$, as
the algorithm of Micali and Vazirani~\cite{MV-80,vaz-90}
(cf.~\cite{blu-90,GT-91}) and solves the b-factor or maximum
degree-constrained subgraph problem in $O(m^{3/2})$ time, similarly
to Gabow~\cite{gab-83}. The logical structure of our algorithm differs
from that of~\cite{MV-80} and sophisticated data structures
(incremental trees for set union~\cite{GT-85}) are used
only in the regular reachability and shortest regular path algorithms
of linear complexity from~\cite{GK-96} (applied as black-box
subroutines) and once in the algorithm for MBP.

Finally, we show that a clique compression technique of Feder and
Motwani~\cite{FM-91} can be extended to certain
skew-symmetric graphs. As a result, our maximum matching algorithm
in a general graph is speeded up to run in
$O(\sqrt{n}m\log(n^2/m)/\log{n})$ time. This matches the best bound for
bipartite matching~\cite{FM-91}.

Fremuth-Paeger and Jungnickel~\cite{FJ-99} developed an algorithm for
MSFP (stated in terms of ``balanced flows'') which combines Dinits'
approach with ideas and tools from~\cite{MV-80,vaz-90}; it runs in
$O(nm^2)$ time for general capacities and in time slightly slower than
$O(\sqrt{n}m)$ in the nonbipartite matching case.

This paper is organized as follows. Basic definitions and facts are
given in Section~\ref{sec:back}. Sections~\ref{sec:theo}
and~\ref{sec:lin} contain theoretical results on
combinatorial and linear programming aspects of IS-flows,
respectively. Section~\ref{sec:gisa} describes Anstee's type algorithm
for MSFP. The shortest regular augmenting path algorithm and a high
level description of the blocking IS-flow method are developed in
Section~\ref{sec:sbf}. Section~\ref{sec:shortp} gives a short review
on facts and algorithms in~\cite{GK-96} for regular path problems.
Using it, Section~\ref{sec:iter} explains the idea of implementation of
a phase in the blocking IS-flow method. It also bounds the number of
phases for special skew-symmetric networks. Section~\ref{sec:acyc}
completes the description of the blocking IS-flow algorithm
by reducing the problem of finding a totally blocking IS-flow in an
acyclic skew-symmetric network to the above-mentioned auxiliary problem
in a usual acyclic digraph and devising a fast procedure to solve the
later. The concluding Section~\ref{sec:mat} discusses implications for
matchings and their generalizations, and explains how to speed up
the implied maximum matching algorithm by use of the clique
compression.

This paper is self-contained up to several quotations from~\cite{GK-96}.
Main results presented in this paper were announced in extended
abstract~\cite{GK-95}. Subsequently, the authors found a flaw in the
original fast implementation of a phase in the blocking IS-flow method.
It was corrected in a revised version of this paper (circulated in 2001)
where problem MBP and its weighted analog were introduced and
efficiently solved,
whereas the original version (identical to preprint~\cite{GK-99})
embraced only the content of Sections~\ref{sec:theo}--\ref{sec:iter}.

\section{\Large Preliminaries}\label{sec:back}

By a {\em skew-symmetric graph} we mean a digraph $G = (V,E)$ with
a mapping (involution) $\sigma$ of $V\cup E$ onto itself
such that:
(i) for each $x\in V\cup E$, $\sigma(x)\ne x$ and $\sigma(\sigma(x))=x$;
(ii) for each $v\in V$, $\sigma(v)\in V$;
and (iii) for each $a=(v,w)\in E$, $\sigma(a)=(\sigma(w),\sigma(v))$.
Although parallel arcs are allowed in $G$,
an arc leaving a node $x$ and entering a node $y$ is denoted by $(x,y)$
when it is not confusing.
We assume that $\sigma$ is fixed (when there are several such mappings)
and explicitly included in the description of $G$.
The node (arc) $\sigma(x)$ is called {\em symmetric} to a node (arc) $x$
(using, for brevity, the term {\em symmetric\/} rather than
skew-symmetric).
Symmetric objects are also called {\em mates}, and we usually use
notation with
primes for mates: $x'$ denotes the mate $\sigma(x)$ of an element
$x$. Note that $G$ can contain an arc $a$ from a node $v$ to its mate
$v'$; then $a'$ is also an arc from $v$ to $v'$.

Unless mentioned otherwise, when talking about paths (cycles), we mean
directed paths (cycles). The symmetry $\sigma$ is extended in a natural
way to paths, subgraphs, and other objects in $G$; e.g., two paths
(cycles) are symmetric if the elements of one of them are symmetric to
those of the other and go in the reverse order. Note that $G$ cannot
contain self-symmetric paths or cycles. Indeed, if
$P=(x_0,a_1,x_1,\ldots,a_k,x_k)$ is such a path (cycle), choose arcs
$a_i$ and $a_j$ such that $i\le j$, $a_j=\sigma(a_i)$ and $j-i$ is
minimum. Then $j>i+1$ (as $j=i$ would imply $\sigma(a_i)=a_i$ and
$j=i+1$ would imply $\sigma(x_i)=x_{j-1}=x_i$). Now
$\sigma(a_{i+1})=a_{j-1}$ contradicts the minimality of $j-i$.

We call a function $h$ on $E$ {\em symmetric\/} if
$h(a) = h(a')$ for all $a \in E$.

A {\em skew-symmetric network\/} is a quadruple $N=(G,\sigma,u,s)$
consisting of a skew-symmetric graph $G=(V,E)$ with symmetry $\sigma$,
a nonnegative integer-valued symmetric function $u$ (of {\em arc
capacities}) on $E$, and a {\em source} $s\in V$.
The mate $s'$ of $s$ is the {\em sink} of $N$. A {\em flow} in $N$
is a function $f: E \to\Rset_+$ satisfying the capacity constraints
$$
    f(a) \leq u(a) \qquad \mbox{for all}\;\; a \in E
$$
and the conservation constraints
$$
   \div_f(x):=\sum_{(x,y) \in E} f(x,y) - \sum_{(y,x) \in E} f(y,x) =0
     \qquad\mbox{for all}\;\; x \in V - \{s,s'\} .
$$
The value $\div_f(s)$ is called the {\em value} of $f$ and denoted by
$|f|$; we usually assume that $|f|\ge 0$. Now {\em IS-flow} abbreviates
{\em integer symmetric flow}, the main object that we study in this
paper. The {\em maximum IS-flow problem (MSFP)} is to find
an IS-flow of maximum value in $N$.

The integrality requirement is important: if we do not require
$f$ to be integral, then for any integer flow $f$ in $N$, the flow $f'$,
defined by $f'(a):=(f(a)+f(a'))/2$ for $a\in E$, is a flow of the same
value as $f$, which is symmetric but not necessarily integral.
Therefore, the {\em fractional} skew-symmetric flow problem is
equivalent to the ordinary flow problem.

Note that, given a digraph $D=(V(D),A(D))$ with two specified nodes
$p$ and $q$ and nonnegative integer capacities of the arcs, we can
construct a skew-symmetric graph $G$ by taking a disjoint copy $D'$
of $D$ with all arcs reversed, adding two extra nodes $s$ and $s'$,
and adding four arcs $(s,p),(s,q'),(q,s'),(p',s')$ of infinite
capacity, where $p',q'$ are the copies of $p,q$ in $D'$,
respectively. Then there is a natural one-to-one correspondence
between integer flows from $p$ to $q$ in $D$
and the IS-flows from $s$
to $s'$ in $G$. This shows that MSFP generalizes the classical
(integer) max-flow problem.

\medskip
\noindent
{\bf Remark.} Sometimes it is useful to consider a sharper version of
MSFP in which double-sided capacity constraints $\ell(a)\le f(a)\le
u(a)$, $a\in E$, are imposed, where $\ell,u:E\to\Zset_+$ and $\ell
\le u$ ({\em problem DMSFP}). Similarly to the max-flow problem
with upper and lower capacities~\cite{FF-62}, DMSFP is reduced
to MSFP
in the skew-symmetric network $N'$ obtained from $N$ by subdividing
each arc $a=(x,y)$ into three arcs $(x,v),(v,w),(w,y)$
with (upper) capacities $u(a),u(a)-\ell(a),u(a)$, respectively, and
adding extra arcs $(s,w)$ and $(v,s')$ with capacity $\ell(a)$ each.
It is not difficult to show (e.g., using Theorem~\ref{tm:m-m})
that DMSFP has a solution if and only if all extra arcs are saturated
by a maximum IS-flow $f'$ for $N'$, and in this case $f'$ induces a
maximum IS-flow for $N$ in a natural way. For details,
see~\cite{FJ-99}.

\medskip
In our study of IS-flows we rely on results for regular paths
in skew-symmetric graphs. A {\em regular path}, or an {\em r-path},
is a path in $G$ that does not contain a pair of symmetric arcs.
Similarly, an {\em r-cycle} is a cycle that does not contain a pair
of symmetric arcs. The {\em r-reachability problem (RP)\/} is to find
an r-path from $s$ to $s'$ or a proof that there is none.
Given a symmetric function of
{\em arc lengths}, the {\em shortest r-path problem (SRP)\/} is to find
a minimum length r-path from $s$ to $s'$ or a proof that there is none.

\begin{figure}[tb]
  \unitlength=1mm
 \begin{center}
  \begin{picture}(120,80)
\put(50,5){\circle{2}}
\put(52,6){$s$}
\put(50,75){\circle{2}}
\put(52,72){$s'$}
\put(30,40){\circle*{1.5}}
\put(40,40){\circle*{1.5}}
\put(50,40){\circle*{1.5}}
\put(0,20){\line(5,-2){50}}
\put(100,20){\line(-5,-2){50}}
\put(0,60){\line(5,2){50}}
\put(100,60){\line(-5,2){50}}
\put(0,20){\line(1,0){100}}
\put(0,60){\line(1,0){100}}
\put(37,8){$A$}
\put(37,69){$A'$}
\put(15,40){\oval(15,25)}
\put(65,40){\oval(15,25)}
\put(90,30){\line(1,0){25}}
\put(90,50){\line(1,0){25}}
\put(90,30){\line(0,1){20}}
\put(115,30){\line(0,1){20}}
\put(13,39){$X_1$}
\put(63,39){$X_k$}
\put(101,39){$M$}
\put(10,34){\vector(1,0){10}}
\put(20,46){\vector(-1,0){10}}
\put(62,34){\vector(1,0){8}}
\put(70,46){\vector(-1,0){8}}
\put(95,35){\vector(0,1){10}}
\put(110,45){\vector(0,-1){10}}
\put(18,32){\vector(2,-3){10}}
\put(28,63){\vector(-2,-3){10}}
\put(60,35){\vector(-1,-2){10}}
\put(50,65){\vector(1,-2){10}}
\put(100,35){\vector(-3,-4){14}}
\put(86,63){\vector(3,-4){14}}
\put(35,65){\vector(0,-1){50}}
\put(45,65){\vector(0,-1){50}}
\put(11,22){$e^1$}
\put(67,22){$e^k$}
\thicklines
\put(15,17){\vector(0,1){13.5}}
\put(15,50){\vector(0,1){13.5}}
\put(65,17){\vector(0,1){13}}
\put(65,50){\vector(0,1){13}}

  \end{picture}
 \end{center}
\caption{ A barrier}
\label{fig:bar}
  \end{figure}
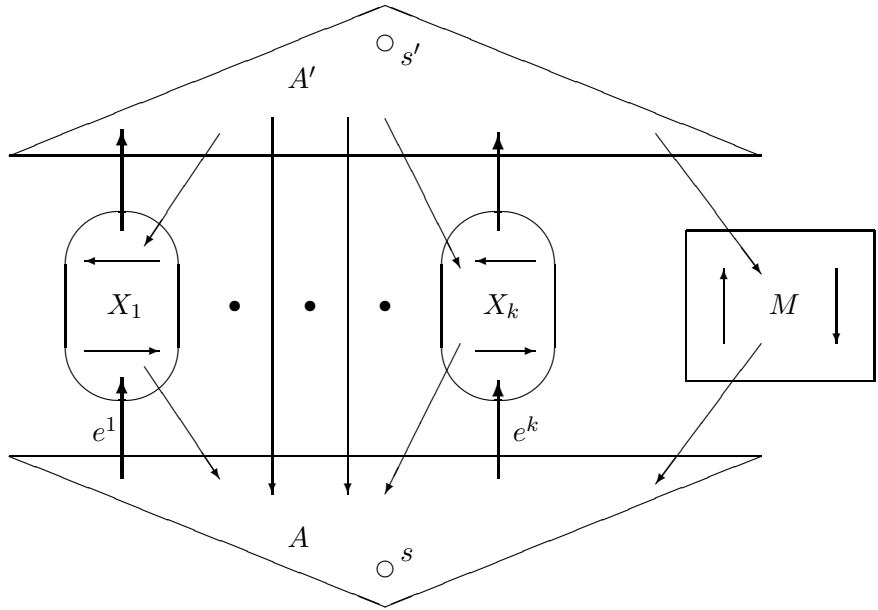

A criterion for the existence of a regular $s$ to $s'$ path is less
trivial than that for the usual path reachability; it involves
so-called barriers. We say that
$$
    \Bscr=(A; X_1, \ldots, X_k)
$$
is an {\em $s$-barrier} if the following conditions hold.
\begin{enumerate}
  \item[(B1)] $A, X_1, \ldots, X_k$ are pairwise disjoint subsets
of $V$, and $s \in A$.
  \spitem[(B2)] For $A' = \sigma(A)$, $A \cap A' = \emptyset$.
  \spitem[(B3)] For $i = 1, \ldots, k$, $X_i$ is self-symmetric, i.e.,
$\sigma(X_i) = X_i$.
  \spitem[(B4)] For $i = 1, \ldots, k$, there is a unique arc, $e^i$,
from $A$ to $X_i$.
  \spitem[(B5)] For $i,j = 1,\ldots, k$ and $i \not = j$, no arc
connects $X_i$ and $X_ j$.
  \spitem[(B6)] For $M := V - (A\cup A' \cup X_1 \cup \ldots \cup X_k)$
and $i = 1, \ldots, k$, no arc connects $X_i$ and $M$.
  \spitem[(B7)] No arc goes from $A$ to $A' \cup M$.
\end{enumerate}
(Note that arcs from $A'$ to $A$, from $X_i$ to $A$, and from $M$ to $A$
are possible.) Figure~\ref{fig:bar} illustrates the definition.
Tutte proved the following (see also~\cite{blu-90,GK-96}).

\begin{theorem}\label{tm:bar}{\rm \cite{tut-67}}
There is an r-path from $s$ to $s'$ if and only if there is no s-barrier.
\end{theorem}

This criterion will be used in Section~\ref{sec:theo}
to obtain an analog of the max-flow min-cut theorem for IS-flows.
RP is efficiently solvable.

\begin{theorem} \label{tm:ratime} {\rm \cite{blu-90,GK-96}}
The r-reachability problem in $G$ can be solved in $O(m)$ time.
\end{theorem}

The methods for the maximum IS-flow problem that we develop apply,
as a subroutine, the r-reachability algorithm of linear
complexity from~\cite{GK-96}, which finds either a regular $s$ to $s'$
path or an $s$-barrier. Another ingredient used in our methods is the
shortest r-path algorithm for the case of nonnegative symmetric
lengths, which runs in $O(m\,\log n)$ time in general, and in $O(m)$
time for all-unit lengths~\cite{GK-96}. The necessary results on RP and
SRP are outlined in Section~\ref{sec:shortp}.

In the rest of this paper, $\sigma$ and $s$
will denote the symmetry map and the source, respectively,
regardless of the network in question, which will allow us to use
the shorter notation $(G,u)$ for a network $(G,\sigma,u,s)$.
Given a simple path $P$, the number of arcs on $P$ is denoted by
$|P|$ and the incidence vector of its arc set in $\Rset^E$
is denoted by $\chi^P$, i.e., $\chi^P(a)=1$ if $a$ is an arc of
$P$, and 0 otherwise.

\subsection{Relationships to Matchings}
\label{ssec:rel_mat}

Given an undirected graph $G' = (V',E')$, a {\em matching} is a subset
$M\subseteq E'$ such that no two edges of $M$ have a common endnode.
The {\em maximum matching problem} is to find a matching $M$
whose cardinality $|M|$ is as large as possible.

There are well-known generalizations of matchings;
for a survey see~\cite{law-76,LP-86,sch-03}. Let
$u_0,u:E'\to\Zset_+\cup\{\infty\}$ and $b_0,b:V'\to\Zset_+$ be
functions such that $b_0\le b$ and $u_0\le u$.
A {\em $(u_0,u)$-capacitated $(b_0,b)$-matching} is a function
$h: E'\to\Zset_+$ satisfying the capacity constraints
$$
    u_0(e)\le h(e)\le u(e) \qquad \mbox{for all}\;\;  e \in E',
$$
and the supply constraints
$$
   b_0(v)\le \sum_{e=\{v,w\}\in E'} h(e) \le b(v) \qquad
        \mbox{for all}\;\;  v \in V'.
$$
The {\em value}\/ of $h$ is defined to be $h(E')$. Hereinafter,
for a numerical function $g$ on a set $S$ and a subset $S'\subseteq
S$, $g(S')$ denotes $\sum_{e\in S'} g(e)$.
Popular special cases are:
a {\em u-capacitated b-matching} (when $b_0=0$);
a {\em degree-constrained subgraph} (when $u\equiv 1$);
a {\em perfect b-matching} (when $u\equiv \infty$ and $b_0=b$);
a {\em b-factor} (when $u\equiv 1$ and $b_0=b$). In these cases one
assigns $u_0=0$.
Typically, in unweighted versions, one is asked for maximizing the
value of $h$ (in the former two cases) or for finding a feasible $h$
(in the latter two cases).

The general maximum $(u_0,u)$-capacitated $(b_0,b)$-matching problem
is reduced to the maximum IS-flow problem (MSFP or DMSFP, depending on
whether both $u_0,b_0$ are zero functions or not) without increasing
the problem size by more than a constant factor. The construction of
the corresponding capacitated skew-symmetric graph $G=(V,E)$
is straightforward (and close to that in~\cite{tut-67}):

(i) for each $v\in V'$, $V$ contains two symmetric nodes
$v_1$ and $v_2$;

(ii) also $V$ contains two additional symmetric nodes
$s$ and $s'$ (the source and the sink);

(iii) for each $e=\{v,w\} \in E'$, $E$ contains two symmetric arcs
$(v_1,w_2)$ and $(w_1, v_2)$ with lower capacity $u_0(e)$ and upper
capacity $u(e)$;

(iv) for each $v\in V'$, $E$ contains two symmetric arcs ($s,v_1)$
and $(v_2,s')$ with lower capacity $b_0(v)$ and upper capacity $b(v)$.

\begin{figure}[tb]
  \unitlength=1mm
 \begin{center}
  \begin{picture}(150,75)
\put(15,15){\circle*{1.5}}
\put(11,14){$a$}
\put(15,30){\circle*{1.5}}
\put(11,29){$b$}
\put(30,45){\circle*{1.5}}
\put(32,44){$c$}
\put(15,60){\circle*{1.5}}
\put(11,59){$d$}
\put(15,75){\circle*{1.5}}
\put(11,74){$e$}
\put(70,45){\circle*{1.5}}
\put(66,44){$s$}
\put(140,45){\circle*{1.5}}
\put(142,44){$s'$}
\put(90,15){\circle*{1.5}}
\put(85,12){$a_1$}
\put(90,30){\circle*{1.5}}
\put(85,27){$b_1$}
\put(90,45){\circle*{1.5}}
\put(86,41){$c_1$}
\put(90,60){\circle*{1.5}}
\put(85,61){$d_1$}
\put(90,75){\circle*{1.5}}
\put(85,75){$e_1$}
\put(120,15){\circle*{1.5}}
\put(122,12){$a_2$}
\put(120,30){\circle*{1.5}}
\put(122,27){$b_2$}
\put(120,45){\circle*{1.5}}
\put(121,41){$c_2$}
\put(120,60){\circle*{1.5}}
\put(122,61){$d_2$}
\put(120,75){\circle*{1.5}}
\put(122,75){$e_2$}
\put(15,30){\line(0,1){30}}
\put(15,60){\line(0,1){15}}
\put(15,30){\line(1,1){15}}
\put(70,45){\line(2,3){20}}
\put(140,45){\line(-2,3){20}}
\put(90,30){\line(2,1){30}}
\put(120,30){\line(-2,1){30}}
\put(90,30){\line(1,1){30}}
\put(120,30){\line(-1,1){30}}
\put(90,60){\line(2,1){30}}
\put(120,60){\line(-2,1){30}}
{                     %thick group
\thicklines
\put(15,15){\line(0,1){15}}
\put(30,45){\line(-1,1){15}}
\put(70,45){\vector(2,-3){19}}
\put(70,45){\vector(4,-3){19}}
\put(70,45){\vector(1,0){19}}
\put(70,45){\vector(4,3){19}}
\put(90,15){\vector(2,1){29}}
\put(90,30){\vector(2,-1){29}}
\put(90,45){\vector(2,1){29}}
\put(90,60){\vector(2,-1){29}}
\put(120,15){\vector(2,3){18.5}}
\put(120,30){\vector(4,3){18.5}}
\put(120,45){\vector(1,0){18.5}}
\put(120,60){\vector(4,-3){18.5}}
}
\put(7,2){Matching}
\put(87,2){Skew-symmetric flow}
   \end{picture}
 \end{center}
\caption{{\sl Reduction example for maximum cardinality
matching}}
\label{fig:red}
  \end{figure}
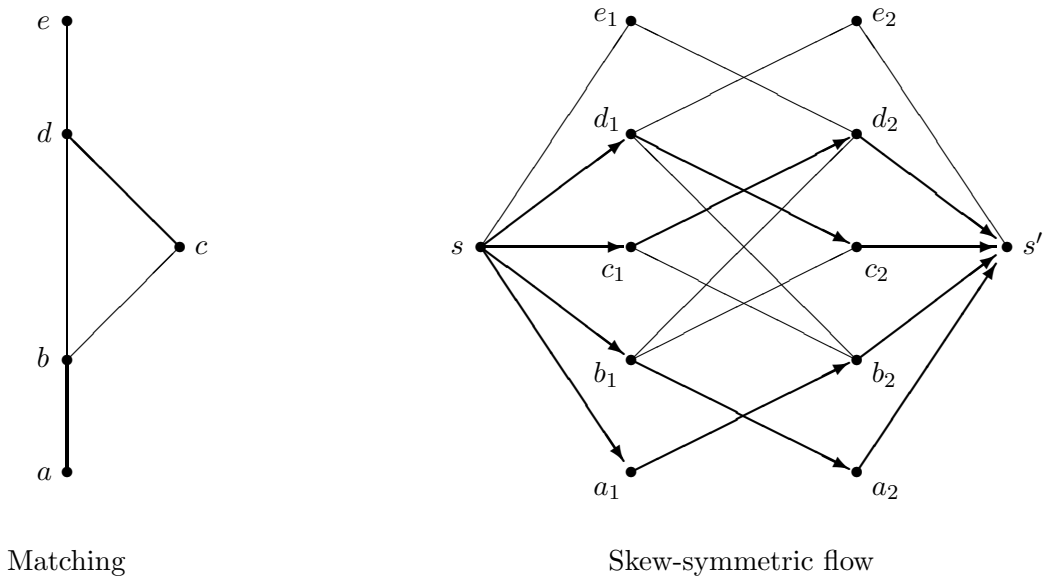

There is a natural one-to-one correspondence between the
$(u_0,u)$-capacitated $(b_0,b)$-matchings $h$ in $G'$ and the IS-flows
$f$ from $s$ to $s'$ in $G$, and the value of $f$ is twice the value
of $h$. Figure~\ref{fig:red} illustrates the correspondence for
matchings.

In case of the b-factor or degree-constrained subgraph problem,
one may assume
that $b$ does not exceed the node degree function of $G$. Therefore,
one can make a further reduction to MSFP in a network with $O(|E'|)$
nodes, $O(|E'|)$ arcs, and unit arc capacities (by getting rid of
lower capacities as in the Remark above and then splitting each arc
$a$ with capacity $q(a)>1$ into $q(a)$ parallel arcs of capacity one).
In Section~\ref{sec:mat} we compare the time bounds of our methods for
MSFP applied to the matching problem and its generalizations above with
known bounds for these problems.

Edmonds and Johnson~\cite{EJ-70} studied the class of integer linear
programs in which the constraint matrix entries are integers
between --2 and +2 and the sum of absolute values of
entries in each column (without including entries from the box
constraints) does not exceed two. Such a problem is often stated in
terms of bidirected graphs (for a survey,
see~\cite[Chapter 36]{sch-03}). Recall that a {\em bidirected graph}
$H=(X,B)$ may contain, besides usual directed edges going from one
node to another, edges directed {\em from} both of its endnodes, and
{\em to} both of them. A particular problem on such an object is the
{\em maximum bidirected flow problem}: given a capacity
function $c:B\to \Zset_+$ and a {\em terminal} $p\in X$,
find a function ({\em biflow}) $g:B\to\Zset_+$
maximizing the value $\div_g(p)$.
(Reasonable versions with more terminals are reduced to this one.)
Here $g\le c$ and $\div_g(x)=0$ for all $x\in X-\{p\}$,
where $\div_g(x)$ is the total biflow on the edges directed
from $x$ minus the total biflow on the edges directed to $x$
(a loop $e$ at $x$ contributes 0, $2g(e)$ or $-2g(e)$).

The maximum IS-flow problem admits a linear time and space reduction
to the maximum biflow problem (in fact, both are equivalent).
More precisely, given an
instance $N=(G=(V,E),\sigma,u,s)$ of MSFP, take a partition $(X,X')$
of $V$ such that $X'=\sigma(X)$ and $s\in X$. For each pair $\{a,a'\}$
of symmetric arcs in $E$ and nodes $x,y\in X$, assign an edge from $x
$ to $y$ if $a$ or $a'$ goes from $x$ to $y$; an edge from both $x,y$
if $a$ or $a'$ goes from $x$ to $\sigma(y)$; an edge to both $x,y$ if
$a$ or $a'$ goes from $\sigma(x)$ to $y$. This produces a bidirected
graph $H=(X,B)$. We set $p:=s$ and assign the capacity $c(e)$ of
each edge $e\in B$ to be the capacity of the arc from which
$e$ is created. There is a one-to-one correspondence between
the IS-flows in $N$ and the
biflows in $(H,c,p)$, and the values of corresponding flows are equal.
A reverse reduction is also obvious. Using these reductions, one can
try to derive results for IS-flows from corresponding results on
biflows, and vice versa.
In this paper we give direct proofs and algorithms for IS-flows.

%----------------------SEC 3
\section{\Large Mini-Theory of Skew-Symmetric Flows}\label{sec:theo}

This section extends the classical flow decomposition, augmenting
path, and max-flow min-cut theorems of Ford and Fulkerson \cite{FF-62}
to the skew-symmetric case.
The {\em support} $\{e\in S: f(e)\ne 0\}$ of a function $f:S\to\Rset$
is denoted by $\supp(f)$.

Let $h$ be a nonnegative integer symmetric function on the arcs of
a skew-symmetric graph $G=(V,E)$.
A path (cycle) $P$ in $G$ is called $h$-{\em regular} if
$h(a)>0$ for all arcs $a$ of $P$ and each arc $a\in P$ such that $a'\in
P$ satisfies $h(a)\ge 2$. Clearly when $h$ is all-unit on $E$, the sets
of regular and $h$-regular paths (cycles) are the same. We call an arc
$a$ of $P$ {\em ordinary} if $a'\not\in P$ and define the $h$-{\em
capacity} $\delta_h(P)$ of $P$ to be the minimum of all values $h(a)$ for
ordinary arcs $a$ on $P$ and all values $\lfloor h(a)/2\rfloor$ for
nonordinary arcs $a$ on $P$.

To state the symmetric flow decomposition theorem, consider an IS-flow
$f$ in a skew-symmetric network $N=(G=(V,E),u)$. An IS-flow $g$ in $N$
is called {\em elementary} if it is representable as $g=\delta\chi^P
+\delta\chi^{P'}$, where $P$ is a simple cycle or a simple path from $s$
to $s'$ or a simple path from $s'$ to $s$, $P'=\sigma(P)$, and
$\delta$ is a {\em positive integer}. Since $g$ is feasible, $P$ is
$u$-regular and $\delta\le\delta_u(P)$. We denote $g$ by
$(P,P',\delta)$. By a {\em symmetric decomposition} of $f$ we mean a set
$D$ of elementary flows such that $f=\sum(g : g\in D)$.
The following {\em symmetric decomposition theorem}
(see~\cite{FJ-99,GK-95}) slightly
generalizes a result by Tutte~\cite{tut-67} that there exists a
symmetric set of $|f|$ paths from $s$ to $s'$ such that any arc $a$
is contained in at most $f(a)$ paths.

\begin{theorem}\label{tm:dec} For an IS-flow $f$ in $G$, there exists a
symmetric decomposition consisting of at most $m$ elementary flows.
\end{theorem}

\begin{proof} We build up an
$f$-regular path $\Gamma$ in $G$ until this path contains a simple
cycle $P$ or a simple path $P$ connecting $s$ and $s'$. This will
determine a member of the desired flow decomposition. Then we
accordingly decrease $f$ and repeat the process for the resulting
IS-flow $f'$, and so on until we obtain the zero flow.

We start with $\Gamma$ formed by a single arc $a\in\supp(f)$. First we
grow $\Gamma$ forward. Let $b=(v,w)$ be the last arc on the current
(simple) path $\Gamma$.
Suppose that $w \not = s,s'$. By the conservation for $f$ at $w$,
$\supp(f)$ must contain an arc $q=(w,z)$. If $q'$ is not on $\Gamma$ or
$f(q)\ge 2$, we add $q$ to $\Gamma$.

Suppose $q'$ is on $\Gamma$ and $f(q)=1$. Let $\Gamma_1$ be the part
of $\Gamma$ between $w'$ and $w$. Then $\Gamma_1$ contains at
least one arc since $w \ne w'$. Suppose there is an arc $\tilde q\in
\supp(f)$ leaving $w$ and different from $q$. Then we can add $\tilde
q$ to $\Gamma$ instead of $q$, forming a longer $f$-regular path.
(Note that since the path $\Gamma$ is simple, ${\tilde q}'$ is not
on $\Gamma$). Now
suppose that such a $\tilde q$ does not exist. Then exactly one unit of
the flow $f$ leaves $w$. Hence, exactly one unit of the flow $f$ enters
$w$, implying that $b=(v,w)$ is the only arc entering $w$ in
$\supp(f)$, and that $f(b)=1$. But $\sigma(d)$ also enters $w$, where
$d$ is the first arc on $\Gamma_1$. The fact that $\sigma(d)\ne b$
(since $\Gamma_1$ is $f$-regular) leads to a contradiction.

Let $(w,z)$ be the arc added to $\Gamma$. If $z$ is not on $\Gamma$,
then $\Gamma$ is a simple $f$-regular path, and we continue growing
$\Gamma$. If $z$ is on $\Gamma$, we discover a simple $f$-regular cycle
$P$.

If $\Gamma$ reaches $s'$ or $s$, we start growing $\Gamma$ backward
from the initial arc $a$ in a way similar to growing it forward. We stop
when an $f$-regular cycle $P$ is found or one of $s$, $s'$ is reached.
In the latter case $P = \Gamma$ is either an $f$-regular path from $s$
to $s'$ or from $s'$ to $s$, or an $f$-regular cycle (containing $s$ or
$s'$).

Form the elementary flow $g=(P,P',\delta)$ with $\delta=\delta_f(P)$ and
reduce $f$ to $f':=f-\delta\chi^P-\delta\chi^{P'}$. Since $P$ is
$f$-regular, $\delta>0$. Moreover, there is a pair $e,e'$ of symmetric
arcs of $P$ such that either $f'(e)=f'(e')=0$ or $f'(e)=f'(e')=1$; we
associate such a pair with $g$. In the former case $e,e'$ vanish in the
support of the new IS-flow $f'$, while in the latter case $e,e'$ can be
used in further iterations of the decomposition process at most once.
Therefore, each pair of arc mates of $G$ is associated with at most two
members of the constructed decomposition $D$, yielding $|D|\le m$.
\end{proof}

The above proof gives a polynomial time algorithm for symmetric
decomposition. Moreover, the above decomposition process can be easily
implemented in $O(nm)$ time, which matches the complexity of standard
decomposition algorithms for usual flows.

The decomposition theorem and the fact that the network has no
self-symmetric cycles imply the following useful property noticed by
Tutte as well.

\begin{corollary}\label{cor:even} {\rm \cite{tut-67}}
For any self-symmetric set $S\subseteq V$
and any IS-flow in $G$, the total flow on the arcs entering $S$, as well
as the total flow on the arcs leaving $S$, is even.
\end{corollary}

\Remark
Another consequence of Theorem~\ref{tm:dec} is that one may assume
that $G$ has no arc entering $s$. Indeed, consider a maximum IS-flow
$f$ in $G$ and a symmetric decomposition $D$ of $f$. Putting together
the elementary flows from $s$ to $s'$ in $D$, we obtain an IS-flow $f'$
in $G$ with $|f'| \geq |f|$, so $f'$ is a maximum flow. Since $f'$ uses
no arc entering $s$ or leaving $s'$, deletion of all such arcs from $G$
produces an equivalent problem in a skew-symmetric graph.

\medskip
Next we state a skew-symmetric version of the augmenting path theorem.
It is convenient to consider the graph $G^+=(V,E^+)$ formed by adding a
reverse arc $(y,x)$ to each arc $(x,y)$ of $G$. For $a \in E^+$, $a^R$
denotes the corresponding reverse arc. The symmetry $\sigma$ is extended
to $E^+$ in a natural way. Given a (nonnegative integer) symmetric
capacity function $u$ on $E$ and an IS-flow $f$ on $G$,
define the {\em residual capacity}
$u_f(a)$ of an arc $a \in E^+$ to be $u(a)-f(a)$ if $a \in E$, and
$f(a^R)$ otherwise. An arc $a \in E^+$ is called {\em residual} if
$u_f(a) > 0$, and {\em saturated} otherwise. Given an IS-flow $g$ in the
network $(G^+,u_f)$, we define the function $f \oplus g$ on $E$ by
setting $(f \oplus g)(a) := f(a) + g(a) - g(a^R)$. Clearly $f \oplus g$
is a feasible IS-flow in $(G,u)$ whose value is $|f|+|g|$.

By an {\em r-augmenting path} for $f$ we mean a $u_f$-regular path from
$s$ to $s'$ in $G^+$. If $P$ is an r-augmenting path and if $\delta \in
\Nset$ does not exceed the $u_f$-capacity of $P$, then we can push
$\delta$ units of flow through a (not necessarily directed) path in $G$
corresponding to $P$ and then $\delta$ units through the path
corresponding to $P'$. Formally, $f$ is transformed into $f \oplus g$,
where $g$ is the elementary flow $(P, P', \delta)$ in $(G^+,u_f)$. Such
an augmentation increases the value of $f$ by $2\delta$.

\begin{theorem}\label{tm:aug} {\rm \cite{tut-67}}
An IS-flow $f$ is maximum if and only if there is no r-augmenting path.
\end{theorem}

\begin{proof}
The direction that the existence of an r-augmenting path implies
that $f$ is not maximum is obvious in light of the above discussion.

To see the other direction, suppose that $f$ is not maximum, and let
$f^*$ be a maximum IS-flow in $G$. For $a \in E$ define $g(a) := f^*(a) -
f(a)$ and $g(a^R) := 0$ if $f^*(a) \geq f(a)$, while $g(a^R) := f(a) -
f^*(a)$ and $g(a) := 0$ if $f^*(a) < f(a)$. One can see that $g$ is a
feasible symmetric flow in $(G^+,u_f)$. Take a symmetric decomposition
$D$ of $g$. Since $|g| = |f^*| - |f| > 0$, $D$ has a
member $(P,P',\delta)$, where $P$ is a $u_f$-regular path from $s$ to
$s'$. Then $P$ is an r-augmenting path for $f$.
\end{proof}

In what follows we will use a simple construction which enables us
to reduce the task of finding an r-augmenting path to the r-reachability
problem. For a skew-symmetric network $(H,h)$, split each arc $a =
(x,y)$ of $H$ into two parallel arcs $a_1$ and $a_2$ from $x$ to $y$
(the {\em first\/} and {\em second split-arcs\/} generated by $a$).
These arcs are endowed with the capacities $[h](a_1) := \lceil h(a)/2
\rceil$ and $[h](a_2) := \lfloor h(a)/2 \rfloor$. Then delete all arcs
with zero capacity $[h]$. The resulting capacitated graph is called the {\em
split-graph} for $(H,h)$ and denoted by $S(H,h)$. The symmetry $\sigma$
is extended to the arcs of $S(H,h)$ in a natural way, by defining
$\sigma(a_i) := (\sigma(a))_i$ for $i=1,2$.

For a path $P$ in $S(H,h)$, its image in $H$ is denoted by $\omega(P)$
(i.e., $\omega(P)$ is obtained by replacing each arc $a_i$ of $P$ by the
original arc $a=:\omega(a_i)$). It is easy to see that if $P$
is regular, then $\omega(P)$ is $h$-regular. Conversely, for any
$h$-regular path $Q$ in $H$, there is a (possibly not unique)
r-path $P$ in $S(H,h)$ such that $\omega(P)=Q$. Indeed, replace
each ordinary arc $a$ of $Q$ by the first split-arc $a_1$ (existing
as $h(a)\ge 1$) and replace each pair $a,a'$ of arc mates in $Q$ by
$a_i,a'_j$ for $\{i,j\}=\{1,2\}$ (taking into account that
$h(a)=h(a')\ge 2$). This gives the required r-path $P$. Thus,
Theorem~\ref{tm:aug} admits the following reformulation in terms
of split-graphs.

\begin{corollary}\label{tm:spaug} An IS-flow $f$ in $(G,u)$ is maximum
if and only if there is no regular
path from $s$ to $s'$ in $S(G^+,u_f)$.
\end{corollary}

Finally, the classic max-flow min-cut theorem states that the maximum
flow value is equal to the minimum cut capacity. A skew-symmetric
version of this theorem involves a more complicated object which is
close to an $s$-barrier occurring in the solvability criterion for the
r-reachability problem given in Theorem~\ref{tm:bar}. We say that
$\Bscr = (A; X_1, \ldots, X_k) $ is an {\em odd $s$-barrier} for
$(G,u)$ if the following conditions hold.

\begin{enumerate}
  \item[(O1)] $A, X_1, \ldots, X_k$ are pairwise disjoint subsets
of $V$, and $s \in A$.
  \spitem[(O2)] For $A' = \sigma(A)$, $A \cap A' = \emptyset$.
  \spitem[(O3)] For $i = 1, \ldots, k$, $X_i$ is self-symmetric, {\it
i.e.,} $\sigma(X_i) = X_i$.
  \spitem[(O4)] For $i = 1, \ldots, k$, the total capacity $u(A, X_i)$
of the arcs from $A$ to $X_i$ is odd.
  \spitem[(O5)] For $i,j = 1, \ldots, k$ and $i \not = j$, no
positive capacity arc connects $X_i$ and $X_j$.
  \spitem[(O6)] For $M := V - (A\cup A' \cup X_1 \cup \ldots \cup X_k)
$ and $ i = 1, \ldots, k$, no positive capacity arc connects $X_i$
and $M$.
  \end{enumerate}
Compare with (B1)--(B7) in Section~\ref{sec:back}.
Define the {\em capacity} $u(\Bscr)$ of $\Bscr$ to be
$u(A, V-A) - k$. Since the source is denoted by $s$ throughout, we
refer to an odd $s$-barrier as {\em odd barrier}.

The following is the {\em maximum IS-flow minimum barrier theorem}.
\begin{theorem}\label{tm:m-m}  {\rm \cite{tut-67}}
The maximum IS-flow value is equal to the minimum odd barrier capacity.
\end{theorem}

\begin{proof}
To see that the capacity of an odd barrier $\Bscr=(A;X_1,\ldots,
X_k)$ is an upper bound on the value of an IS-flow $f$, consider a
symmetric decomposition $D$ of $f$. For each member $g=(P,P',\delta)$ of
$D$, where $P$ is a path from $s$ to $s'$, take the {\em last} arc
$a=(x,y)$ of the {\em first} path $P$ such that $x\in A$. If $y\in
A'$, then the symmetric arc $a'$ (which is in $P'$) also goes from $A$
to $A'$ (by (O2)), and therefore, $g$ uses at least $2\delta$ units of the
capacity of arcs from $A$ to $A'$. Associate $g$ with the pair $a,a'$.
Now let $y\not\in A'$. Since $y\not\in A$, $y$ is either in $Y:=M$ or
in $Y:=X_i$ for some $i$. The choice of $a$ and (O1),(O5),(O6) imply
that $P$ leaves $Y$ by an arc $b$ from $Y$ to $A'$.
Then the symmetric arc $b'$ (which is in $P'$) goes from $A$ to $Y$
(since $Y$ is self-symmetric),
and therefore, $g$ uses at least $2\delta$ units of the capacity
$u(A,Y)$. Associate $g$ with the pair $a,b'$ (possibly $a=b'$).
Note that at least one unit of each capacity $u(A,X_i)$ is not used
under the canonical way we associate the elementary $s$ to $s'$ flows
of $D$ with arcs from $A$ to $V-A$ (since $u(A,X_i)$ is odd, by (O4)).
By these reasonings, $|f|\le u(\Bscr)$.

Next we show that the two values in the theorem are equal. Let $f$ be a
maximum IS-flow. By Corollary~\ref{tm:spaug}, the split-graph
$S=S(G^+,u_f)$ contains no $s$ to $s'$ r-path, so it must contain an
$s$-barrier $\Bscr = (A; X_1, \ldots, X_k)$, by Theorem~\ref{tm:bar}.

Let $e^i$ be the (unique) arc from $A$ to $X_i$ in $S$ (see (B4) in
Section 2). By the construction of $S$, it follows that the residual
capacity $u_f$ of every arc from $A$ to $X_i$ in $G^+$ is zero except
for the arc $\omega(e^i)$, whose residual capacity is one.
Hence,
\begin{itemize}
\item[(i)] if $e^i$ was formed by splitting an arc $a\in E$, then
$a$ goes from $A$ to $X_i$, and $f(a)=u(a)-1$;
\spitem[(ii)] if $e^i$ was formed by splitting $a^R$ for $a\in E$,
then $a$ goes from $X_i$ to $A$, and $f(a)=1$;
\spitem[(iii)] all arcs from $A$ to $X_i$ in $G$, except $a$ in case (i),
are saturated by $f$;
\spitem[(iv)] all arcs from $X_i$ to $A$ in $G$, except $a$ in case (ii),
are free of flow.
\end{itemize}

Furthermore, comparing arcs in $S$ and $G$, we observe that:
\begin{itemize}
\item[(v)] property (B7) implies that the arcs from $A$
to $A'\cup M$ are saturated and the arcs from $A'\cup M$ to $A$ are free
of flow;
\spitem[(vi)] property (B5) implies (O5) and (B6) implies (O6).
\end{itemize}

Properties (i)--(iv),(O5),(O6) together with Corollary~\ref{cor:even}
provide (O4). So $\Bscr$ is an odd $s$-barrier in $G$. We have
$|f|=f(A,V-A)-f(V-A,A)=u(A,V-A)-k$ (in view of (i)--(v)). Hence,
$|f|=u(\Bscr)$.
   \end{proof}

%-----------------------SECTION 4
\section{\Large Integer and Linear Programming
Formulations}\label{sec:lin}

Although methods of solving MSFP developed in subsequent sections will
not use explicitly linear programming aspects exhibited in this
section, such aspects help to understand more about the structure of
IS-flows.

MSFP is stated as an integer program in a straightforward way.
We use function rather than vector notation. For functions $g,h$ on a
set $S$, $g\cdot h$ denotes the inner product $\sum_{x\in S}g(x)h(x)$.
Assuming that no arc of $G$ enters the source $s$ (see the Remark in
the previous section), MSFP can be written as follows:
\begin{eqnarray}
 \mbox{\bf maximize}\;\; |f| = \sum\nolimits_{(s,v) \in E} f(s,v)
   & &  \mbox{\bf subject to} \label{eq:1} \\
 f(a) \geq 0 & & \forall a \in E \label{eq:2} \\
 f(a) \leq u(a) & &  \forall a \in E \label{eq:3} \\
 -\sum\nolimits_{(u,v) \in E} f(u,v)
           + \sum\nolimits_{(v,w) \in E} f(v,w) = 0
   & &   \forall v \in V - \{s,s'\} \label{eq:4} \\
 f(a) - f(\sigma(a)) = 0 & & \forall a \in E \label{eq:5} \\
 f(a) \;\;\;\mbox{integer} & & \forall a \in E \label{eq:6}
\end{eqnarray}

A linear programming formulation for MSFP is obtained by
replacing the integrality condition (\ref{eq:6}) by linear
constraints related to certain objects that we call odd fragments
in $G$. The correctness of the resulting
linear program will be shown by use of the max-min relation
between IS-flow and odd barriers in Theorem~\ref{tm:m-m}.
Alternatively, one can try to derive it from a linear programming
characterization of integer bidirected flows in~\cite{EJ-70} (using the
reduction as in Section~\ref{sec:back}).

An {\em odd fragment} is a pair $\rho = (V_\rho, U_\rho)$, where
$V_\rho$ is a {\em self-symmetric} set of nodes with $s\not\in V_\rho$,
and $U_\rho$ is a subset of arcs entering $V_\rho$ such that the total
capacity $u(U_\rho)$ is odd. The {\em characteristic function}
$\chi_\rho$ of $\rho$ is the function on $E$ defined by

\begin{equation}\label{eq:ch_of}
   \chi_\rho (a) := \left\{
    \begin{array}{rl}
  1 & \mbox{if} \;\; a \in U_\rho \cup \sigma(U_\rho),\\
  -1 & \mbox{if} \;\; a \in \delta(V_\rho) - (U_\rho \cup
       \sigma(U_\rho)),\\
  0 & \mbox{otherwise}.
  \end{array}
        \right.
\end{equation}
Here $\delta(V_\rho)$ is the set of arcs with one end in $V_\rho$ and
the other in $V - V_\rho$. We denote the set of odd fragments by $\Omega$.

Let $f$ be a (feasible) IS-flow and $\rho\in\Omega$. By~(\ref{eq:ch_of})
and the symmetry of $u$, we have $f\cdot\chi_\rho\le
u(U_\rho)+u(\sigma(U_\rho))=2u(U_\rho)$.
Moreover, $f\cdot\chi_\rho$ is at most $2u(U_\rho)-2$;
this immediately
follows from Corollary~\ref{cor:even} and the fact
that $u(U_\rho)$ is odd. This gives new linear constraints
for MSFP:
\begin{equation}
  f\cdot \chi_\rho\le 2u(U_\rho)-2\quad \mbox{for each}\;\; \rho\in\Omega.
   \label{eq:8}
  \end{equation}
Addition of these constraints enables us to drop off the symmetry
constraints (\ref{eq:5}) and the integrality constraints (\ref{eq:6})
without changing the optimum value of the linear program. This fact is
implied by the following theorem.

\begin{theorem}\label{tm:opt}
Every maximum IS-flow is an optimal solution to the linear program
{\rm (\ref{eq:1})--(\ref{eq:4}), (\ref{eq:8})}.
\end{theorem}

\begin{proof}
Assign a dual variable $\pi(v)\in\Rset$ (a {\it potential}) to each node
$v\in V$, $\gamma(a)\in\Rset_+$ (a {\it length}) to each arc $a\in E$,
and $\xi(\rho)\in\Rset_+$ to each odd fragment $\rho\in\Omega$.
Consider the linear program:
 \begin{eqnarray}
  \mbox{\bf minimize}\;\;  \psi(\pi,\gamma,\xi):=\sum_E u(a)\gamma(a)+
  \sum_{\Omega}(2u(U_\rho)-2)\xi(\rho)
         & & \mbox{\bf subject to} \label{eq:9} \\
  \gamma(a)\ge 0 & & \forall a\in E \label{eq:10} \\
  \xi(\rho)\ge 0 & & \forall \rho\in\Omega \label{eq:11} \\
  \pi(s)=0 & & \label{eq:12} \\
  \pi(s')=1 & & \label{eq:13} \\
  \pi(v)-\pi(w)+\gamma(a)+\sum_{\Omega}\xi(\rho)\chi_\rho(a)
      \ge 0 & & \forall a=(v,w)\in E .
     \label{eq:14}
\end{eqnarray}

In fact, (\ref{eq:9})--(\ref{eq:14}) is dual to linear
program~(\ref{eq:1})--(\ref{eq:4}),(\ref{eq:8}). (To see this,
introduce an extra arc $(s',s)$, add the conservation
constraints for $s$ and $s'$, and replace the objective (\ref{eq:1})
by $\max\{f(s',s)\}$. The latter generates the dual constraint
$\pi(s')-\pi(s)\ge 1$. We can replace it by the equality or
impose~(\ref{eq:12})--(\ref{eq:13}).) Therefore,
 \begin{equation}
  \max\;|f| =\min\; \psi(\pi,\gamma,\xi), \label{eq:15}
 \end{equation}
where the maximum and minimum range over the corresponding feasible
solutions.

We assert that every maximum IS-flow $f$ achieves the maximum in
(\ref{eq:15}). To see this, choose an odd barrier $ \Bscr = (A; X_1,
\ldots, X_k) $ of minimum capacity $u(\Bscr)$. For $i=1,\ldots,k$,
let $U_i$ be the set of arcs from $A$ to $X_i$; then $\rho_i=(X_i,U_i)$
is an odd fragment for $G,u$. Define $\pi(v)$ to be 0 for $v\in A$, 1
for $v\in A'$, and $1/2$ otherwise. Define $\gamma(a)$ to be 1 for $a\in
(A,A')$, $1/2$ for $a\in (A,M)\cup (M,A')$, and 0 otherwise, where
$M=V-(A\cup A'\cup X_1\cup\ldots\cup X_k)$. Define
$\xi(\rho_i)=1/2$ for $i=1,\ldots,k$, and $\xi(\rho)=0$ for the other
odd fragments in $(G,u)$.

One can check that (\ref{eq:14}) holds for all arcs $a$ (e.g.,
both values $\pi(w)-\pi(v)$ and $\gamma(a)+
\sum_{\Omega}\xi(\rho)\chi_\rho(a)$ are equal to 1
for $a=(v,w)\in (A,A')$, and 1/2 for $a=(v,w)\in (A,L)\cup (L,A')$,
where $L:=V-(A\cup A')$). Thus $\pi,\gamma,\xi$ are feasible.

Using the fact that $u(A,M)=u(M,A')$, we observe that
$u\cdot\gamma=u(A,A')+u(A,M)$. Also
  $$
   \sum_{\Omega}(2u(U_\rho)-2)\xi(\rho)= \sum_{i=1}^k
   \frac{1}{2}(2u(U_{i})-2)=\left(\sum_{i=1}^ku(A,X_k)\right)-k.
  $$
This implies $\psi(\pi,\gamma,\xi)=u(\Bscr)$, and now the result
follows from Theorem~\ref{tm:m-m}.
\end{proof}

%----------------------------SEC 5
\section{\Large Algorithm Using a Good Pre-Solution}
\label{sec:gisa}

Anstee~\cite{ans-85,ans-87} developed efficient methods for b-factor
and b-matching problems (unweighted or weighted) based on the idea
that a good pre-solution can easily be found by solving a
corresponding flow problem.
In this section we adapt his approach to solve
the maximum IS-flow problem in a skew-symmetric network $N=
(G=(V,E),u)$. The algorithm that we devise is relatively simple; it
finds a ``nearly optimal'' IS-flow and then makes $O(n)$ augmentations
to obtain a maximum IS-flow. The algorithm consists of four stages.

The {\em first\/} stage ignores the fact that $N$ is skew-symmetric
and finds an integer maximum flow $g$ in $N$ by use of a
max-flow algorithm. Then we set $h(a):=(g(a)+g(a'))/2$
for all arcs $a\in E$. Since $\div_h(s)=\div_g(s)/2-\div_g(s')/2
=\div_g(s)$, $h$ is a maximum flow as well. Also
$h$ is symmetric and {\em half-integer}.
Let $Z$ be the set of arcs on which $h$ is not integer. If
$Z=\emptyset$, then $h$ is already a maximum IS-flow; so assume
this is not the case.

The {\em second\/} stage applies simple transformations of $h$ to
reduce $Z$.
Let $H=(X,Z)$ be the subgraph of $G$ induced by $Z$. Obviously,
for each $x\in V$, $\div_h(x)$ is an integer,
so $x$ is incident to an even number of
arcs in $Z$. Therefore, we can decompose $H$ into simple, not
necessarily directed, cycles $C_1,\ldots,C_r$ which are pairwise
arc-disjoint. Moreover, we can find, in linear time, a decomposition
in which each cycle $C_i$ is either self-symmetric ($C_i=\sigma(C_i)$)
or symmetric to another cycle $C_j$ ($C_i=\sigma(C_j)$).

To do this, we start with some node $v_0\in X$ and grow in $H$ a
simple (undirected) path $P=(v_0,e_1,v_1,\ldots,e_q,v_q)$ such that
the mate $v'_i$ of each node $v_i$ is not in $P$. At each step, we
choose in $H$ an arc $e\ne e_q$ incident to the last node $v_q$
($e$ exists since $H$ is eulerian); let $x$ be the other end node of
$e$. If none of $x,x'$ is in $P$, then we add $e$ to $P$. If some
of $x,x'$ is a node of $P$, $v_i$ say, then we shorten $P$ by
removing its end part from $e_{i+1}$ and delete from $H$ the arcs
$e_{i+1},\ldots,e_q,e$ and their mates. One can see that the arcs
deleted induce a self-symmetric cycle (when $x'=v_i$) or two
disjoint symmetric cycles (when $x=v_i$). We also remove the
isolated nodes created by the arc deletions and change the
initial node $v_0$ if needed.
Repeating the process for the new current graph $H$ and path $P$, we
eventually obtain the desired decomposition $\Cscr$, in $O(|Z|)$
time.

Next we examine the cycles in $\Cscr$. Each pair $C,C'$ of
symmetric cycles is canceled by sending a half unit of flow
through $C$ and through $C'$, i.e., we increase (resp. decrease)
$h(e)$ by 1/2 on each forward (resp.
backward) arc $e$ of these cycles. The resulting function $h$ is
symmetric, and $\div_h(x)$ is preserved at each node $x$, whence $h$
is again a maximum symmetric flow. Now suppose that two
self-symmetric cycles $C$ and $D$ meet at a node $x$. Then they
meet at $x'$ as well. Concatenating the $x$ to $x'$ path in $C$ and
the $x'$ to $x$ path in $D$ and concatenating the rests of $C$ and $D
$, we obtain a pair of symmetric cycles and cancel these cycles as
above.
These cancellations result in $\Cscr$ consisting of pairwise
{\em node-disjoint} self-symmetric cycles, say $C_1,\ldots,C_k$.
The second stage takes $O(m)$ time.

The {\it third\/} stage transforms $h$ into an IS-flow $f$ whose
value $|f|$ is at most $k$ units below $|h|$. For each
$i$, fix a node $t_i$ in $C_i$ and change $h$ on $C_i$ by sending a
half unit of flow through the $t_i$ to $t'_i$ path in $C_i$ and
through the reverse to the $t'_i$ to $t_i$ path in it. The resulting
function $h$ is integer and symmetric and the divergences preserve
at all nodes except for the nodes $t_i$ and $t'_i$
where we have $\div_h(t_i)=-\div_h(t'_i)=1$ for each $i$
(assuming, without loss of generality, that all $t_i$'s are
different from $s'$). Therefore, $h$ is, in essence, a
multiterminal IS-flow with sources $s,t_1,\ldots,t_k$ and
sinks $s',t'_1,\ldots,t'_k$. A genuine IS-flow $f$ from $s$
to $s'$ is extracted by reducing $h$ on some $h$-regular paths.
More precisely, we add to $G$ artificial arcs $e_i=(s,t_i)$, $i=1,
\ldots, k$ and their mates, extend $h$ by ones to these
arcs and construct a symmetric decomposition $\Dscr$ (defined in
Section~ \ref{sec:theo}) for the obtained function $h'$ in the
resulting graph $G'$ (clearly $h'$ is an IS-flow of value
$|h|+k$).

Let $\Dscr'$ be the set of elementary flows in $\Dscr$ formed
by the paths or cycles which contain artificial arcs. Then
$\delta=1$ for each $(P,P',\delta)\in\Dscr'$. Define $f':=h'-
\sum(\chi^P+\chi^{P'}: (P,P',1)\in\Dscr')$. Then $f'$ is an
IS-flow in $G'$, and $|f'|\ge|h'|-2k\ge|h|-k$.
Moreover, since $f(e_i)=0$ for
$i=1,\ldots,k$, the restriction $f$ of $f'$ to $E$ is an IS-flow in
$G$, and $|f|=|f'|$. Thus, $|f|\ge|h|-k$, and now the facts that
$k\le n/2$ (as the nodes $t_1,\ldots,t_k,t'_1,\ldots,t'_k$ are
different) and that $h$ is a maximum flow in $N$ imply that the
value of $f$ differs from the maximum IS-flow value by $O(n)$.
The third stage takes $O(nm)$ time (the time
needed to construct a symmetric decomposition of $h'$).

The final, {\em fourth}, stage transforms $f$ into a maximum
IS-flow. Each iteration applies the r-reachability algorithm (RA)
mentioned in Section~\ref{sec:back} to the split-graph
$S(G^+,u_f)$ in order to find a $u_f$-regular $s$ to $s'$ path
$P$ in $G^+$ and then augment the current IS-flow $f$ by the
elementary flow
$(P,P',\delta_{u_f}(P))$ as explained in Section~\ref{sec:theo}.
Thus, a maximum IS-flow in $N$ is constructed in $O(n)$ iterations.
Since the RA runs in $O(m)$ time (by Theorem~\ref{tm:ratime}), the
fourth stage takes $O(nm)$ time.

Summing up the above arguments, we conclude with the following.
  \begin{theorem} \label{tm:gisa}
  The above algorithm finds a maximum IS-flow in $N$ in
$O(M(n,m)+nm)$ time, where $M(n,m)$ is the running time of the
max-flow procedure it applies.
  \end{theorem}

%------------------------SEC 6
\section{\Large Shortest R-Augmenting Paths and Blocking IS-Flows}
\label{sec:sbf}

Theorem \ref{tm:aug} and Corollary \ref{tm:spaug} prompt an
alternative method
for finding a maximum IS-flow in a skew-symmetric network
$N=(G,u)$, which is analogous to the method of Ford and Fulkerson for
usual flows. It starts with the zero flow, and at each iteration, the
current IS-flow $f$ is augmented by an elementary flow in $(G^+,u_f)$
(found by applying the r-reachability algorithm to $S(G^+,u_f)$).
Since each iteration increases the value of $f$ by at least two,
a maximum IS-flow is constructed in pseudo-polynomial time. In
general, this method is not competitive to the method of
Section~\ref{sec:gisa}.

More efficient methods involve the concepts of
shortest r-augmenting paths and shortest blocking IS-flows that
we now introduce. Let $g$ be an IS-flow in a
skew-symmetric network $(H=(V,W),h)$. We call $g(W)$
($=\sum_{e\in W} g(e)$) the {\em volume\/} of $g$. Considering a
symmetric decomposition $D=\{(P_i,P'_i,\delta_i): i=1,\ldots,k)$ of
$g$, we have
 $$
   g(W)=\sum(\delta_i|P_i|+\delta_i|P'_i| : i=1,\ldots,k) \ge
    |g|\min\{|P_i|: i=1,\ldots,k\}.
$$
This implies
  \begin{equation}
   g(W)\ge |g|\mbox{\rdist}_{S(H,h)}(s,s'),    \label{eq:5-1}
  \end{equation}
where $\rdist_{H'}(x,y)$ denotes the minimum length of a regular
$x$ to $y$ path in a skew-symmetric graph $H'$ (the {\em regular
distance\/} from $x$ to $y$). We say that an IS-flow $g$ is
\begin{itemize}
\item[(i)] {\em shortest\/} if (\ref{eq:5-1}) holds with equality,
i.e., some (equivalently, any) symmetric decomposition of $g$
consists of shortest $h$-regular paths from $s$ to $s'$;

\spitem[(ii)] {\em totally blocking\/} if there is no $(h-g)$-regular
path from $s$ to $s'$ in $H$, i.e., we cannot augment $g$ using only
residual capacities in $H$ itself;

\spitem[(iii)] {\em shortest blocking\/} if $g$ is shortest (as in (i))
and
\begin{equation}
   \rdist_{S(H,h-g)}(s,s') > \rdist_{S(H,h)}(s,s'). \label{eq:5-2}
\end{equation}
\end{itemize}

Note that a shortest blocking IS-flow is not necessarily totally
blocking, and vice versa.

Given a skew-symmetric network $N=(G,u)$, the {\em shortest
r-augmenting path method (SAPM)\/}, analogous to the method of
Edmonds and Karp~\cite{EK-72} for usual flows, starts with the zero
flow, and each iteration augments the current IS-flow $f$ by a
shortest elementary flow $g=(P,P',\delta_{u_f}(P))$.

The {\em shortest blocking IS-flow method (SBFM)\/}, analogous to
Dinits' method~\cite{din-70}, starts with the zero flow, and each
{\em phase} (big iteration) augments the current IS-flow
$f$ by performing the following two steps:
  \begin{itemize}
   \item[(P1)] Find a shortest blocking IS-flow $g$ in $(G^+,u_f)$.
   \item[(P2)] Update $f:=f\oplus g$.
  \end{itemize}

Both methods terminate when $f$ no longer admits r-augmenting paths
(i.e., $g$ becomes the zero flow). The following observation is
crucial for our methods.

\begin{lemma}\label{lm:incr}
Let $g$ be a shortest IS-flow in $(G^+,u_f)$, and let $f':=f\oplus g$.
Let $k$ and $k'$ be the minimum lengths of r-augmenting paths for $f$
and $f'$, respectively. Then $k'\ge k$. Moreover, if $g$ is a shortest
blocking IS-flow, then $k'>k$.
\end{lemma}

\begin{proof}
Take a shortest $u_{f'}$-regular path $P$ from $s$ to $s'$ in $G^+$.
Then $|P|=k'$ and $g'=(P,P',1)$ is an elementary flow in $(G^+,u_{f'})$.

Note that $\supp(g)$ does not contain opposed arcs $a=(x,y)$ and
$b=(y,x)$. Otherwise decreasing $g$ by one on each of $a,b,a',b'$
(which are, obviously, distinct), we would
obtain the IS-flow $\tilde g$ in $(G^+,u_f)$ such that $|\tilde g|=|g|$
and $\tilde g(E^+)<g(E^+)$, which is impossible because $\tilde
g(E^+)\ge k|\tilde g|$ and $g(E^+)=k|g|$. This implies that each arc
$a$ in the set $Z:=\{a\in E^+: g(a^R)=0\}$ satisfies
  \begin{equation}
     u_{f'}(a)=u_f(a)-g(a). \label{eq:5-3}
  \end{equation}

If $\supp(g')\subseteq Z$, then $g'$ is a feasible IS-flow in
$(G^+,u_f)$ (by~(\ref{eq:5-3})), whence $k'=g'(E^+)/|g'|\ge k$.
Moreover, if, in addition, $g$ is a shortest blocking IS-flow, then
(\ref{eq:5-2}) and the fact that $g'\le u_f-g$ (by (\ref{eq:5-3}))
imply $k'>k$.

Now suppose there is an arc $e\in E^+$ such that $g'(e)>0$ and
$g(e^R)>0$. For each $a\in E^+$, put $\lambda(a):=\max\{0,g(a)+g'(a)
-g(a^R)-g'(a^R)\}$. One can check that $\lambda(a)\le u_f(a)$ for all
arcs $a$ and that $\div_\lambda(v)=0$ for all nodes $v\ne s,s'$.
Therefore, $\lambda$ is an IS-flow in $(G^+,u_f)$ with
$|\lambda|=|g|+|g'|=|g|+2$. Also $\lambda(E^+)<g(E^+)+g'(E^+)$
since for the $e$ above, $\lambda(e)+\lambda(e^R)<g'(e)+g(e^R)$. We
have
  $$
   2k'=g'(E^+)>\lambda(E^+)-g(E^+)\ge k(|g|+2)-k|g|=2k,
  $$
yielding $k'>k$.
\end{proof}

Thus, each iteration of SAPM does not decrease the minimum
length of an r-augmenting path, and each phase of
SBFM increases this length. This gives upper bounds on the
numbers of iterations.
  \begin{corollary}\label{cor:sapm}
SAPM terminates in at most $(n-1)m$ iterations.
  \end{corollary}

(This follows by observing, in the proof of Lemma~\ref{lm:incr},
that on the iterations with the same length of shortest
r-augmenting paths, the subgraph of $G^+$ induced by the arcs
contained in such paths is monotone nonincreasing, and each
iteration reduces the capacity of some arc of this subgraph, as well
as the capacity of its mate, to zero or one.)

  \begin{corollary}\label{cor:n-1}
SBFM terminates in at most $n-1$ phases.
  \end{corollary}

As mentioned above, SBFM can be considered as a skew-symmetric
analog of Dinits' blocking flow algorithm. Recall that each phase
of that algorithm constructs a blocking flow in the
subnetwork $H$ formed by the nodes and arcs of shortest augmenting
paths. Such a network is acyclic (moreover, layered), and a
blocking flow in $H$ is easily constructed in $O(nm)$ time.

The  problem of finding a shortest blocking IS-flow ((P1) above) is more
complicated. Let $H$ be the subgraph of $G^+$ formed by the nodes and
arcs contained in shortest $u_f$-regular $s$ to $s'$ paths. Such an $H$
need not be acyclic (a counterexample is not difficult).
In Section~\ref{sec:iter} we will show
that problem (P1) can be reduced to a seemingly easier task, namely,
to finding a totally
blocking IS-flow in a certain acyclic network $(\bar H,\bar h)$.
Such a network arises when the shortest r-path algorithm
from~\cite{GK-96} is applied to the split-graph $S(G^+,u_f)$ with
unit arc lengths.
First, however, we need to tell more about the r-reachability and
shortest r-path algorithms from~\cite{GK-96}.

%---------------------------SEC 7
\section{\Large Properties of Regular and Shortest Regular Path
Algorithms}
\label{sec:shortp}

In this section we exhibit certain properties of the algorithms
in~\cite{GK-96}, referring the reader to that paper for details. We
also establish an additional fact (Lemma~\ref{lm:acyc}), which will be
used later.

%------------SUBSEC 1
\subsection{The Regular Reachability Algorithm (RA)}\label{sec:ra}

Let $\Gamma=(V,E)$ be a skew-symmetric graph with source $s$ and sink
$s'=\sigma(s)$ (as before, $\sigma$ is the symmetry map).
A {\em fragment} (or an $s$-{\em fragment}) in $\Gamma$ is a pair
$\phi=(V_\phi,e_\phi = (v,w))$, where $V_\phi $ is a self-symmetric
set of nodes of $\Gamma$ with $s\not\in V_\phi$ and $e_\phi$ is an arc
entering $V_\phi$, i.e., $v\not\in V_\phi\ni w$ (cf. the definition of
odd fragments in Section~\ref{sec:lin}).
We refer to $e_\phi$ and $e'_\phi$ as the {\em base} and
{\em anti-base} arcs of $\phi$, respectively. Let us
say that the fragment is {\em well-reachable} if

(i) for each node $x\in V_\phi$, there is an r-path from
$w$ to $x$ in the subgraph induced by $V_\phi$ (and therefore, an
r-path from $x$ to $w'= \sigma(w)$), and

(ii) there is an r-path from $s$ to $v$ disjoint from $V_\tau$.

\begin{figure}[tb]
  \unitlength=1mm
 \begin{center}
  \begin{picture}(150,70)
\multiput(0,0)(80,0){2}{
\put(15,15){\circle*{1.5}}
\put(13,17){$z$}
\put(35,10){\circle*{1.5}}
\put(37,9){$v$}
\put(35,25){\circle*{1.5}}
\put(37,22){$w$}
\put(5,35){\circle*{1.5}}
\put(4,31.5){$x$}
\put(65,35){\circle*{1.5}}
\put(64,31){$y'$}
\put(5,45){\circle*{1.5}}
\put(4,47){$y$}
\put(65,45){\circle*{1.5}}
\put(64,47){$x'$}
\put(35,55){\circle*{1.5}}
\put(30,56){$w'$}
\put(35,70){\circle*{1.5}}
\put(30.5,68){$v'$}
\put(55,65){\circle*{1.5}}
\put(56,62){$z'$}
{\thicklines
\put(35,10){\vector(0,1){14}}
\put(35,55){\vector(0,1){14}}
}
}                                % end of multiput
\put(25,35){\circle*{1.5}}
\put(27,36){$c$}
\put(45,30){\circle*{1.5}}
\put(42,31){$b'$}
\put(25,50){\circle*{1.5}}
\put(27,47){$b$}
\put(45,45){\circle*{1.5}}
\put(41,41){$c'$}
\put(33,38){$V_\phi$}
\put(35,40){\oval(30,40)}
\put(36.5,15){$e_\phi$}
\put(36.5,64){$e'_\phi$}
\put(15,15){\vector(2,1){18.5}}
\put(45,30){\vector(4,1){18.5}}
\put(25,35){\vector(-1,0){18.5}}
\put(5,45){\vector(4,1){18.5}}
\put(65,45){\vector(-1,0){18.5}}
\put(35,55){\vector(2,1){18.5}}
\put(35,25){\vector(-1,1){9.0}}
\put(45,45){\vector(-1,1){9.0}}
\put(95,15){\vector(1,2){19.5}}
\put(115,25){\vector(1,2){19.5}}
\put(115,25){\vector(-3,1){28.5}}
\put(115,25){\vector(3,1){28.5}}
\put(85,45){\vector(3,1){28.5}}
\put(145,45){\vector(-3,1){28.5}}
%\qbesier(115,25)(105,40)(115,55)
\put(113.5,30){\vector(0,1){20}}
\put(116.5,30){\vector(0,1){20}}
\put(20,0){Before trimming}
\put(100,0){After trimming}
   \end{picture}
 \end{center}
\caption{Fragment trimming example}
\label{fig:trim}
  \end{figure}
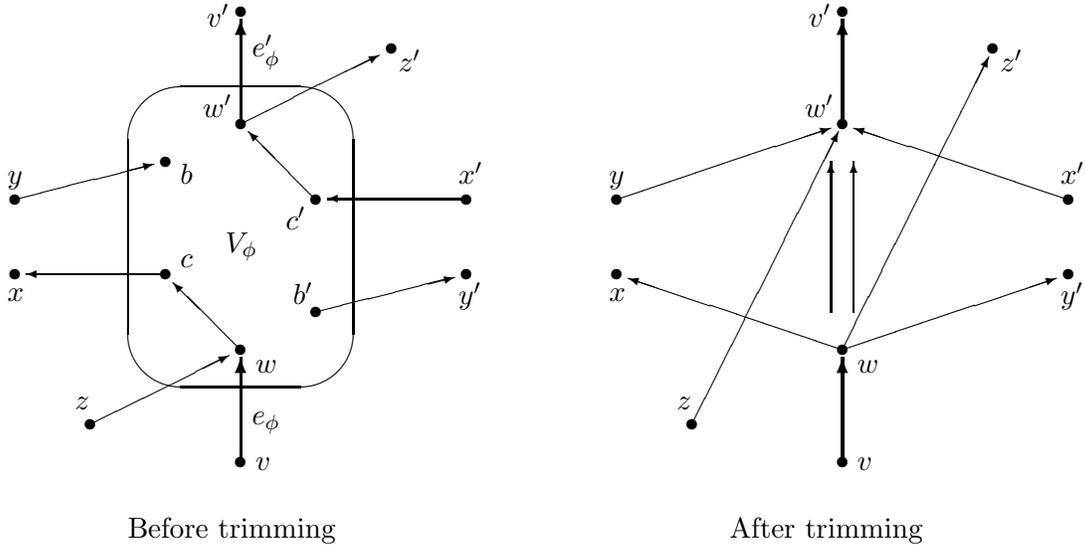

The {\em trimming operation} applied to $\phi$ (which is analogous to
shrinking a blossom in matching algorithms) transforms $\Gamma$ by
removing the nodes of $V_\phi -\{w, w'\}$ and modifying the arcs as
follows.
  \begin{itemize}
    \item[(T1)]
Each arc $a=(x,y) \in E$ such that either
$x,y \in V - V_\phi$ or $a = e_\phi$ or $a=e'_\phi$
remains an arc from $x$ to $y$.
    \spitem[(T2)]
Each arc $(x,y) \in E-\{e'_\phi\}$ that leaves $V_\phi$ is replaced
by an arc from $w$ to $y$, and each arc $(x,y) \in E-\{e_\phi\}$ that
enters $V_\phi$ is replaced by an arc from $x$ to $w'$.
    \spitem[(T3)]
Each arc with both ends in $V_\phi$ is replaced
by an arc from $w$ to $w'$.
  \end{itemize}
(A variant of trimming deletes all arcs in (T3).) The image of an
arc $a$ in the new graph is denoted again by $a$ (so its end nodes can
be changed, but not its name). Figure~\ref{fig:trim} illustrates
fragment trimming. The new $\Gamma$ is again skew-symmetric.

The algorithm RA relies on the following property.

\begin{statement}\label{st:sym} {\rm \cite{GK-96}}
If $\phi$ is a well-reachable fragment, then trimming $\phi$ preserves
the existence (or non-existence) of a regular path from $s$ to $s'$.
\end{statement}

RA searches for a regular $s$ to $s'$ path in $\Gamma$, starting with
the trivial path $P=s$. Each iteration
either increases the current r-path $P$, or reveals a well-reachable
fragment and trims it, producing the new current graph $\Gamma$ and
accordingly updating $P$. It terminates when either an $s$ to $s'$
r-path $\bar P$ or an $s$-barrier $\bar\Bscr$ in the final graph
$\bar\Gamma$ is found (cf. Theorem~\ref{tm:bar}). The
postprocessing stage extends $\bar P$ into a regular $s$ to $s'$ path
$P$ of the initial $\Gamma$ (cf. Statement~\ref{st:sym}) in the
former case (the {\em path restoration procedure}) and extends
$\bar\Bscr$ into a barrier $\Bscr$ of the initial $\Gamma$ in the
latter case (the {\em barrier restoration procedure}).

The fragments of current graphs revealed by RA determine fragments of
the initial $\Gamma$ in a natural way; all fragments are well-reachable.
Moreover, the set $\Phi$ of these fragments of the initial $\Gamma$ is
{\em well-nested}. This means that
  \begin{itemize}
   \item[(F1)] for distinct $\phi,\psi\in\Phi$, either
$V_\phi\subset V_\psi$ or $V_\psi\subset V_\phi$ or
$V_\phi\cap V_\psi=\emptyset$, and
   \spitem[(F2)] for $\phi,\psi\in\Phi$, if $V_{\psi}\subset V_\phi
$ and $e_{\psi}\in\delta(V_\phi)$ then $e_{\psi}=e_\phi$, and
if $V_{\phi}\cap V_{\psi}=\emptyset$ and $e_{\psi}\in\delta(V_{\phi})$
then $e_{\phi}\not\in\delta(V_{\psi})$.
  \end{itemize}
(Recall that for $X\subseteq V$, $\delta(X)$ is the set of arcs with
one end in $X$ and the other in $V-X$.) Let us say that a path $R$
in $\Gamma$ is {\em compatible with} $\Phi$ if for each $\phi\in\Phi$,
$p:=|R\cap\delta(V_\phi)|\le 2$, and if $p=2$ then $R$ contains
exactly one of $e_\phi,e'_\phi$. The following additional properties
(relying on (T1)--(T3)) are important:
  \begin{myitem}
any regular $s$ to $s'$ path $\bar P$ in the final graph $\bar\Gamma$
is extendable to a regular $s$ to $s'$ path $P$ compatible with $\Phi$
in the initial $\Gamma$, and the path restoration procedure applied to
$\bar P$ constructs such a $P$ in $O(|P|+d)$ time, where $d$ is the total
size of maximal fragments in $\Phi$ traversed by $P$;
  \label{eq:P}
  \end{myitem}
\vspace{-5pt}
  \begin{myitem}
for each $\phi\in\Phi$ and each arc $a\ne e'_\phi$ leaving $V_\phi$,
there exists a compatible with $\Phi$ r-path $Q_\phi(a)$ with the first
arc $e_\phi$, the last arc $a$ and all intermediate nodes in $V_\phi$;
such a path can be constructed by (a phase of) the path
restoration procedure in $O(|V_\phi|)$ time.
  \label{eq:Q}
  \end{myitem}

A fast implementation of RA (supported by the disjoint set union data
structure of~\cite{GT-85}) runs in linear time, as indicated in
Theorem~\ref{tm:ratime}.

%-----------------SUBSEC 2
\subsection{The Shortest Regular Path Algorithm
(SRA)}\label{sec:srpa}

We now consider the shortest regular path problem (SRP) in a
skew-symmetric graph $\Gamma=(V,E)$ with {\em nonnegative
symmetric} lengths $\ell(e)$ of the arcs $e\in E$: find a minimum
length regular path from $s$ to $s'$. One may assume that $s'$ is
r-reachable from $s$. The dual problem involves above-mentioned
fragments. Define the characteristic function $\chi_\phi$ of a fragment
$\phi = (V_\phi, e_\phi)$ by
 \begin{equation}\label{eq:ch_f}
  \chi_\phi(a):= \left\{
  \begin{array}{rl}
  1  & \mbox{for $a=e_\phi, e'_\phi$}, \\
  -1 & \mbox{for $a \in \delta(V_\phi) - \{e_\phi, e'_\phi\} $}, \\
  0  & \mbox{for the remaining arcs of $\Gamma$}.
  \end{array}
   \right.
 \end{equation}
(Compare with (\ref{eq:ch_of}).)
For a function $\pi: V\to\Rset$ (of node {\em potentials})
and a nonnegative function $\xi$ on a set $\Phi$ of fragments,
define the {\em reduced length} of an arc $e = (x,y)$ to be
$$
  \lxipi (e): = \ell(e) + \pi(x) - \pi(y) +
  \sum_{\phi \in \Phi} \xi(\phi)\chi_\phi(e) .
$$

An optimality criterion for SRP can be formulated as follows.
  \begin{theorem} {\rm \cite{GK-96}}\label{tm:sp}
A regular path $P$ from $s$ to $s'$ is a shortest r-path if and
only if there exist a potential $\pi: V\to\Rset$,
a set $\Phi$ of fragments, and a {\em positive} function $\xi$ on
$\Phi$ such that
\begin{eqnarray}
  \lxipi (e) \geq 0 & \mbox{for each} & e \in E; \label{eq:6-2} \\
  \lxipi (e) = 0 & \mbox{for each} & e\in P;  \label{eq:6-3}\\
  \chi^P \cdot \chi_\phi = 0 & \mbox{for each} & \phi\in\Phi.
      \label{eq:6-4}
\end{eqnarray}
\end{theorem}

The {\em shortest r-path algorithm (SRA)}
from \cite{GK-96} implicitly maintains $\pi,\Phi,\xi$ in the input
graph $\Gamma$ and iteratively modifies the graph by trimming certain
fragments.
Let $\Gamma^0$ be the subgraph of the current $\Gamma$ with the same
set of nodes and with the arcs having zero reduced length, called the
current {\em 0-subgraph} (recall that the arcs of the current graph
are identified with the corresponding arcs of the initial one). Each
iteration applies the above r-reachability algorithm RA to search for a
regular $s$ to $s'$ path in $\Gamma^0$. If such a path is found, the
algorithm terminates and outputs this path to a postprocessing stage.
If such a path does not exist, then, using the
barrier $\Bscr=(A;X_1,\ldots,X_k)$ in $\Gamma^0$ constructed by RA,
the iteration trims the fragments determined by the sets $X_i$ and
updates $\pi,\Phi,\xi$, modifying $\Gamma^0$. The reduced lengths
of the arcs within the newly and previously extracted fragments, as
well as of their base and anti-base arcs, are not changed.

Let $\bar\Gamma$ and $\bar\Gamma^0$ denote the final graph $\Gamma$ and
the 0-subgraph in it, respectively, and $\bar P$ the regular $s$ to $s'$
path in $\bar\Gamma^0$ found by the algorithm.
Let $\Gamma^0$ stand for the 0-subgraph of the initial graph $\Gamma$
(concerning the reduced arc lengths determined by the resulting
$\pi,\Phi,\xi$). We call $\Gamma^0$ and $\bar\Gamma^0$ the {\em full}
and {\em trimmed 0-graphs}, respectively.
The postprocessing stage applies the path restoration
procedure of RA to extend $\bar P$ into a regular $s$ to $s'$ path $P$
in $\Gamma^0$, in time indicated in~(\ref{eq:P}).
It also explicitly constructs $\Gamma^0$ (in linear time).

Note that any $s$ to $s'$ r-path or r-cycle $Q$ in $\Gamma$ compatible
with $\Phi$ satisfies $\chi^Q\cdot\chi_\phi=0$ for each $\phi\in\Phi$.
By~(\ref{eq:P}), {\em any} $s$ to $s'$ r-path $\bar P$ in
$\bar\Gamma^0$ is extendable to an $s$ to $s'$ r-path $P$ in
$\Gamma^0$ compatible with $\Phi$. Therefore, $P$ is shortest, by
Theorem~\ref{tm:sp}.

\begin{theorem} \label{tm:sharp} {\rm \cite{GK-96}}
  For nonnegative symmetric arc lengths $\ell$,
SRA runs in $O(m \log n)$ time, and in $O(m \sqrt{\log L})$ time
when $\ell$ is integer-valued and $L$ is the maximum arc length.
Furthermore, the algorithm constructs (implicitly) $\pi,\Phi,\xi$ as
in Theorem~\ref{tm:sp}, where $\pi$ is {\em anti-symmetric} (i.e.,
$\pi(x)=-\pi(x')$ for all $x\in V$), and constructs (explicitly) the
trimmed 0-graph $\bar\Gamma^0$ and the full 0-graph $\Gamma^0$
such that:
  \begin{itemize}
   \spitem[{\rm (A1)}] $\Phi$ is well-nested (obeys
{\rm (F1)--(F2)}) and consists of well-reachable fragments in
$\Gamma^0$; in particular, $\lxipi(e_\phi)=0$ for each $\phi\in\Phi$;
   \spitem[{\rm (A2)}] $\Phi$ satisfies~(\ref{eq:P}) and~(\ref{eq:Q})
with $\Gamma^0,\bar\Gamma^0$ instead of $\Gamma,\bar\Gamma$;
in particular, any regular $s$ to $s'$ path of $\bar\Gamma^0$ is
(efficiently) extendable to a shortest regular $s$ to $s'$ path in
$(\Gamma,\ell)$.
  \end{itemize}
\end{theorem}

(Note that the anti-symmetry of $\pi$ and the symmetry of $\ell$ and
$\chi_\phi$ for all $\phi\in\Phi$ imply that the reduced length
function $\lxipi$ is symmetric. Therefore, the graphs $\Gamma^0$ and
$\bar\Gamma^0$ are indeed skew-symmetric.)
Let $\Phi^{\max}$ denote the set of maximal fragments in $\Phi$.
The sets $V_\phi$ for $\phi\in\Phi^{\max}$
are pairwise disjoint (by (F1)), and the graph $\bar\Gamma^0$ can be
directly obtained from $\Gamma^0$ by simultaneously trimming the
fragments in $\Phi^{\max}$.

In the next section we will take advantage of the relationship
between r-paths in $\bar\Gamma^0$ and shortest r-paths in
$(\Gamma,\ell)$ indicated in (A2). Another important property of
$\bar\Gamma^0$ is as follows.

\begin{lemma} \label{lm:acyc}
If the length $\ell(C)$ of every cycle $C$ in $\Gamma$ is
positive, then $\bar\Gamma^0$ is acyclic. In particular,
$\bar\Gamma^0$ is acyclic if all arc lengths are positive.
\end{lemma}

\begin{proof}
Suppose $\bar\Gamma^0$ contains a (not necessarily regular) simple
cycle $\bar C$. In view of~(\ref{eq:Q}), $\bar C$ is extendable to a
cycle $C$ of $\Gamma^0$ compatible with $\Phi$. Then
$\chi^C\cdot\chi_\phi=0$ for all $\phi\in\Phi$. This implies that
the original length $\ell(C)$ and the reduced length $\lxipi(C)$
are the same (since the changes in $\lxipi$ due to $\pi$ cancel out
as we go around the cycle). Since all arcs of $C$ have zero reduced
length, $\ell(C)=\lxipi(C)=0$. This contradicts the hypotheses of the
lemma.
\end{proof}

%--------------------------- SEC 8
\section{\Large Reduction to an Acyclic Network and Special Cases}
\label{sec:iter}

We continue the description of the shortest blocking IS-flow method
(SBFM) for solving the maximum IS-flow problem in a
network $N=(G=(V,E),u)$ begun in Section~\ref{sec:sbf}. Let $f$ be
a current IS-flow in $N$. We show that the task of finding a shortest
blocking IS-flow $g$ in $(G^+,u_f)$ (step (P1) of a phase
of SBFM) reduces to finding a totally blocking IS-flow
in an acyclic network.

Build the split-graph $\Gamma=S(G^+,u_f)$ and apply the above
shortest regular path algorithm to $\Gamma$ with the
{\em all-unit} length function $\ell$ on the arcs. It constructs
$\phi,\Phi,\xi$ as in Theorems~\ref{tm:sp} and \ref{tm:sharp},
taking $O(m)$ time (since $L=1$). SRA also constructs the trimmed
0-graph $\bar\Gamma^0$, the main object we will deal with.
By Lemma~\ref{lm:acyc}, $\bar\Gamma^0$ is
acyclic. Also the following property takes place.

\begin{lemma}\label{lm:ub}
  Let $a\in E^+$ be an arc with $u_f(a)>1$, and let $a_1,a_2$ be
the corresponding split-arcs in $\Gamma$. Then
$\lxipi(a_1)=\lxipi(a_2)$. Moreover, none of $a_1,a_2$ can be the
base or anti-base arc of any fragment in $\Phi$.
  \end{lemma}
  \begin{proof}
Since $a_1,a_2$ are parallel arcs, for each $\phi\in\Phi$, $a_1$
enters (resp. leaves) $V_\phi$ if and only if $a_2$ enters (resp.
leaves) $V_\phi$. This implies that $\lxipi(a_1)\ne\lxipi(a_2)$ can
happen only if one of $a_1,a_2$ is the base or anti-base arc of
some fragment in $\Phi$. Suppose $a_1\in\{e_\phi,e'_\phi\}$ for
some $\phi\in\Phi$ (the case $a_2\in\{e_\phi,e'_\phi\}$ is similar).
Then $\lxipi(a_1)=0$ (by (A1) in Theorem~\ref{tm:sharp}). Using
property~(F2) from Section~\ref{sec:shortp} (valid as $\Phi$ is
well-nested), one can see that $a_2$ is not the base or anti-base
arc of any fragment in $\Phi$. Therefore,
$\chi_{\psi}(a_2)\le\chi_{\psi}(a_1)$ for all $\psi\in\Phi$,
yielding $\lxipi(a_2)\le\lxipi(a_1)$. Moreover, the latter
inequality is strict because $\chi_\phi(a_2)=-1<1=\chi_\phi(a_1)$
and $\xi(\phi)>0$. Now $\lxipi(a_1)=0$ implies $\lxipi(a_2)<0$,
contradicting (\ref{eq:6-2}).
  \end{proof}

Let $E^0\subseteq E^+$ be the set of (images of) zero
reduced length arcs of $\Gamma$. Lemma~\ref{lm:ub} implies that
the base arc $e_\phi$ of each fragment $\phi\in\Phi$ in $\Gamma$
is generated by an arc $e\in E^0$ with $u_f(e)=1$. We can
identify these $e$ and $e_\phi$ and consider $\phi$ as a fragment
of $G^+$ as well. One can see that $\bar
\Gamma^0$ is precisely the split-graph for $(\bar H,\bar h)$, where
$\bar H=(\bar V,\bar E^0)$ is obtained from $H=(V,E^0)$ by trimming
the maximal fragments in $\Phi$, and $\bar h$ is the restriction of
$u_f$ to $\bar E^0$.

Based on the property of each fragment to have unit capacity of
the base arc, we reduce step (P1) to the desired problem, namely:
  \begin{itemize}
  \item[(B)] {\em Find a totally blocking IS-flow in $(\bar H,\bar
h)$}.
  \end{itemize}

To explain the reduction,
suppose we have found a solution $\bar g$ to (B).
For each maximal fragment $\phi$ in $\Phi$ with
$e_\phi\in\supp(\bar g)$, we
have $\bar g(e_\phi)=1$; therefore, exactly one unit of flow goes
out of the head of $e_\phi$, through an arc $a\in \bar E^0$ say.
We choose the path $Q=Q_\phi(a)$ as in~(\ref{eq:Q}) to connect
$e_\phi$ and $a$ in (the subgraph on $V_\phi$ of) $H$ and then
push a unit of flow
through $Q$ and a unit of flow through the symmetric path $Q'$.
Doing so for all maximal fragments $\phi$, we extend $\bar g$ to
an IS-flow $g$ in $(H,h)$, where $h$ is the restriction of $u_f$ to
$E^0$. Moreover, $g$ is a shortest blocking IS-flow in $(G^+,u_f)$.

Indeed, the fact that the chosen paths $Q$ have zero reduced length and
are compatible with $\Phi$ implies that a symmetric decomposition of
$g$ consists of shortest $u_f$-regular paths (cf. (A2) in
Theorem~\ref{tm:sharp}); so $g$ is shortest. Also $G^+$ cannot contain
a $(u_f-g)$-regular $s$ to $s'$ path $R$ of length $g(E^+)/|g|$. For
such an $R$ would be a path in $H$ compatible with $\Phi$ (in view of
Theorem~\ref{tm:sp}); then the arcs of $R$ occurring in $\bar H$ should
form an $(\bar h-\bar g)$-regular $s$ to $s'$ path in it, contrary to
the fact that $\bar g$ is totally blocking.

Since each path $Q_\phi(a)$ is constructed in $O(|V_\phi|)$ time
(by~(\ref{eq:Q})), and the sets $V_\phi$ of maximal fragments $\phi$ are
pairwise disjoint, the reduction to (B) takes linear time.

  \begin{lemma}\label{lm:redu}
A totally blocking IS-flow in $(\bar H,\bar h)$ can be extended to
a shortest blocking IS-flow in $(G^+,u_f)$, in  $O(m)$ time.
\qedt
  \end{lemma}
  \begin{corollary}\label{cor:cvb}
SBFM solves the maximum IS-flow problem in $O(qT(n,m)+qm)$ time,
where $q$ is the number of phases ($q\le n$) and $T(n,m)$
is the time needed to find a totally blocking IS-flow in an acyclic
network with at most $n$ nodes and $m$ arcs.
  \end{corollary}

Clearly $T(n,m)$ is $O(m^2)$, as a totally blocking flow can be
constructed by $O(m)$ applications of the regular reachability
algorithm; this is slower compared with the phase time $O(nm)$ in
Dinits' algorithm. However, we shall show in the next section that
problem (B) can be solved in $O(nm)$ time as well. Moreover,
the bound will be better for important special cases.

Next we estimate the number of phases.
For the standard max-flow problem,
the number of phases of Dinits' algorithm becomes significantly less
than $n$ in the cases of unit arc capacities and unit ``node
capacities''. To combine these into one case, given a network
$N=(G=(V,E),u)$ with integer capacities $u$, for a node $x\in V$,
define the {\em transit capacity}\/ $u(x)$ to be the minimum of values
$\sum_{y:(x,y)\in E} u(x,y)$ and $\sum_{y:(y,x)\in E}u(y,x)$.
Define
 $$
   \Delta:=\Delta(N):=\sum(u(x) : x\in V-\{s,s'\}).
 $$
As shown in \cite{ET-75,kar-73-2}), the number $q$ of
phases of the blocking flow method does not exceed $2\sqrt{\Delta}$.
In particular, if $u\equiv 1$ then $q=O(\sqrt{m})$, and if the transit
capacities $u(x)$ of all nodes $x\ne s,s'$ ({\em inner} nodes)
are ones, e.g., in the case arising from the bipartite matching
problem, then $q=O(\sqrt{n})$.

A similar argument works for skew-symmetric networks (see
also~\cite{FJ-99} for a special case).
   \begin{lemma}\label{lm:root}
The number of phases of SBFM is at most
$\min\{n,2\sqrt{\Delta}\}$.
   \end{lemma}
   \begin{proof}
After performing $d:=\sqrt{\Delta}$ phases, the r-distance
from $s$ to $s'$ in the network $N'=(G^+,u_f)$ for the current IS-flow
$f$ becomes greater than $d$, by Lemma~\ref{lm:incr}. Let $f^*$
be a maximum IS-flow in $N$, and let $g$ be defined as in the proof of
Theorem~\ref{tm:aug}. Then $g$ is a feasible IS-flow in $N'$ and
$|g|=|f^*|-|f|$. We assert that $|g|\le d$, which immediately implies
that the number of remaining phases is at most $d/2$, thus
proving the lemma. To see this, take a symmetric decomposition
$\Dscr$ of $g$ consisting of elementary flows $(P,P',\delta)$
with $\delta=1$. Let $\Dscr'$ be the family of $s$ to $s'$ paths
$P,P'$ in $\Dscr$; then $|\Dscr'|\ge |g|$. It is easy to see that
$u_f(x)=u(x)$ for each inner node $x$. Each path
in $\Dscr'$ contains at least $d$ inner nodes, and therefore, it
uses at least $d$ units of the total transit capacity of inner nodes of
$N'$. So we have $d|\Dscr'|\le \Delta(N)$. This implies $|g|\le d$.
   \end{proof}

%-----------------------SEC 9
\section{\Large Finding a Totally Blocking IS-Flow in an Acyclic
          Network}
   \label{sec:acyc}

Our aim is to show the following.
 \begin{theorem} \label{tm:ac-unit}
A totally blocking IS-flow in an acyclic skew-symmetric graph with
$O(n)$ nodes, $O(m)$ arcs, and unit arc capacities can be found in
$O(n+m)$ time.
  \end{theorem}

This together with Corollary~\ref{cor:cvb} and Lemma~\ref{lm:root}
yields the following result for the shortest blocking IS-flow method.
\begin{corollary} \label{cor:unit}
In case of a network $N$ with unit arc capacities, SBFM can be
implemented so that it finds a maximum IS-flow in
$O(m\sqrt{\Delta(N)})$ time (assuming $n=O(m)$). In
particular, if the indegree or outdegree of each node is at most one,
then the running time becomes $O(\sqrt{n}m)$.
  \end{corollary}

In the second half of this section we will extend
Theorem~\ref{tm:ac-unit} to general capacities, in which case the phase
time will turn into $O(nm)$, similarly to Dinits' algorithm.

For convenience we keep the original notation for the network in
question. Let $G=(V,E)$ be a skew-symmetric {\em acyclic} graph with
source $s$ and the capacity $u(e)=1$ of each arc $e\in E$.
One may assume that each node belongs to a path from $s$ to
$\sigma(s)$.

First of all we make a reduction to the {\em maximal balanced path-set
problem (MBP)} stated in the Introduction. Since
$G$ is acyclic, one can assign, in linear time, a potential function
$\pi:V\to\Zset$ which is {\em antisymmetric} ($\pi(x)=-\pi(\sigma
(x))$ for each $x\in V$) and {\em increasing} on the arcs
($\pi(y)>\pi(x)$ for each $(x,y)\in E$).
(Indeed, a function $q:V\to\Zset$ increasing on the arcs is constructed,
in linear time, by use of the standard topological sorting.
Now set $\pi(v):=q(v)-q(\sigma(v))$, $v\in V$.)
Subdivide each arc $(x,y)$ with $\pi(x)<0$ and $\pi(y)>0$ into two arcs
$(x,z)$ and $(z,y)$ and assign zero potential to $z$.
The new graph $G$, with $O(m)$ nodes and $O(m)$ arcs, is again
skew-symmetric, and the problem remains essentially the same.

Let $\Gamma$ be the subgraph of the new $G$ induced by the nodes
with nonnegative potentials. Then $\Gamma\cup\sigma(\Gamma)=G$ and
$\Gamma\cap\sigma(\Gamma)=(Z,\emptyset)$, where $Z$ is the
self-symmetric set of zero potential nodes of $G$. Also $\Gamma$
contains $\sigma(s)$.

Clearly every $s$ to $\sigma(s)$ path $P$ of $G$ meets $Z$ at exactly one
node $z$, which subdivides $P$ into an $s$ to $z$ path $R'$ in
$\sigma(\Gamma)$ and a $z$ to $\sigma(s)$ path $Q$ in $\Gamma$.
Then $P$ is regular if and only if $\sigma(R')$ and
$Q$ are arc-disjoint. Conversely, let $Q,R$ be two arc-disjoint $Z$ to
$\sigma(s)$ paths in $\Gamma$ beginning at symmetric nodes in $Z$. Then
the concatenation $\sigma(Q)\cdot R$ (as well as $\sigma(R)\cdot Q$) is
a regular $s$ to $\sigma(s)$ path of $G$ if and only if $Q$ and $R$ are
arc-disjoint.

This shows that our particular totally blocking IS-flow problem is
reduced, in linear time, to MBP with $\Gamma,\sigma(s),Z$ (in fact, the
problems are equivalent).
Theorem~\ref{tm:ac-unit} is implied by the following.
  \begin{theorem} \label{tm:bbp}
MBP is solvable in linear time.
  \end{theorem}

We devise an algorithm for MBP and prove Theorem~\ref{tm:bbp}.
Let the input of MBP consist of an acyclic graph
$\Gamma=(X,U)$, a sink $t$ and a source set $Z$ with a map
(involution) $\sigma:Z\to Z$ giving a partition of $Z$ into pairs.
We say that two arc-disjoint $Z$ to $t$ paths $Q,R$ beginning at
``symmetric'' sources $z,\sigma(z)$ form a {\em good pair}, and say
that a collection of pairwise arc-disjoint $Z$ to $s'$ paths in
$\Gamma$ is a {\em balanced path-set} if its members can be partitioned
into good pairs. So the task is to find a maximal (or ``blocking'')
balanced path-set.

Each iteration of the algorithm will reduce the
arc set of $\Gamma$ and, possibly, the set $Z$, and we sometimes will
use index $i$ for objects in the input of $i$-th iteration. So
$\Gamma_1=(X_1,U_1)$ is the initial graph. Without loss of generality,
one may assume that initially each source has zero indegree and
  \begin{itemize}
\item[(C1)] each node of $\Gamma$ lies on a path from $Z$ to $t$,
  \end{itemize}
and will maintain these properties during the algorithm.

The iteration input will include a path $D$ from a certain node of
$\Gamma$ to $t$, called the {\em pre-path}.
Initially, $D$ is trivial: $D=t$.
The nodes of $\Gamma$ not in $Z\cup\{t\}$ are called {\em inner}.
The current $\Gamma$ may contain special inner nodes, called
{\em complex} ones. They arise when the algorithm shrinks certain
subgraphs of $\Gamma$; the initial graph has no complex nodes.
The adjacency structure of $\Gamma$ is given by double-linked
lists $I_x$ and $O_x$ of the incoming and outgoing arcs,
respectively, for each node $x$.
The arc set of a path $P$ is denoted by $E(P)$.

An $i$-th iteration begins with extending $D$ to a
$Z$ to $t$ path $P$ in a natural way;
this takes $O(|P|-|D|)$ time.
Let $z$ be the first node of $P$. Then we try to obtain a good pair
by constructing a path from $z'=\sigma(z)$ to $t$, possibly
rearranging $P$. By standard arguments, a good pair for $z,z'$
exists if and only if there exists a path $A$ from $z'$ to $t$,
with possible backward arcs, in which the forward arcs belong to
$U-E(P)$ and the backward arcs belong to $E(P)$, called an {\em
augmenting path} w.r.t. $P$. For certain reasons, we admit $A$
to be self-intersecting in nodes (but not in arcs). Once $A$ is found,
the symmetric difference $E(P)\triangle E(A)$ gives a good pair $Q,R$
(taking into account that $\Gamma$ is acyclic).

To search for an augmenting path, we replace each arc $e=(x,y)\in
E(P)$ by the reverse arc $\bar e=(y,x)$; let $\bar \Gamma=(X,\bar U)$
be the resulting graph, and $\bar P$ the $t$ to $Z$ path
reverse to $P$. Thus, we have to construct a (directed) path from $z'$
to $t$ in $\bar\Gamma$ or establish that it does not exist.

To achieve the desired time bound, we apply a variant of
depth first search which we call here {\em transit depth first search
(TDFS)} (such a search procedure was applied in~\cite{kar-70}).
The difference from the standard depth first search (DFS) is as
follows. When scanning a new outgoing arc $(x,y)$ in the list $O_x$
of a current node $x$, if $y$ has already been visited, then DFS
stays at $x$. In contrast, TDFS moves from $x$ to $y$, making $y$ the
new current node. Both procedures maintain the stack
of arcs traversed only in forward direction and ordered by the
time of their traversal. If all outgoing arcs of the current node
$x$ are already scanned, then the last arc $(w,x)$ of the stack is
traversed in backward direction and $w$ becomes the new current
node. We refer to the path determined by the stack, from the initial
node to the current one, as the {\em active path}. Note that in case of
TDFS the active path may be self-intersecting (while it is simple in
DFS).

We impose the condition that the outgoing arc lists of
$\bar\Gamma$ are arranged so that
  \begin{itemize}
\item[(C2)] for each node $x\ne z$ of $\bar P$, the arc $\bar e$
of $\bar P$ leaving $x$ is the {\em last} element of $\bar O_x$.
  \end{itemize}
This guarantees that TDFS would scan $\bar e$ after all other outgoing
arcs of $x$ (i.e., the arcs of $\bar P$ are ignored as long as
possible).

At an iteration, we apply TDFS to $\bar\Gamma$ starting from $z'$ as
above. The search terminates
when either it reaches $t$ or it returns to $z'$ having all arcs of
$O_{z'}$ traversed. In the first case ({\em breakthrough}) the final
active path $\bar A$ in $\bar\Gamma$ determines the desired augmenting
path $A$ in $\Gamma$, and we create a good pair $Q,R$ as described
above. In the second case, the non-existence of a good pair for the
given $z,z'$ is declared. Consider both cases.

\medskip
{\em Breakthrough case.} Delete from $\Gamma$ the arcs of $Q,R$ and
then delete all the nodes and arcs that are no longer contained in $Z$
to $t$ paths (thus maintaining (C1)). This is carried out by an obvious
{\em cleaning procedure} in $O(q)$ time, where $q$ is the number of
arcs deleted. If $Q$ or $R$ contains a
complex node, the iteration finishes by transforming $Q,R$ into a good
pair of paths of the initial graph; this is carried out by a
{\em path expansion procedure} which will be described later.
The obtained $\Gamma,Z$ form the input of the next iteration, and
the new pre-path $D$ is assigned to be the trivial path $t$.
If $\Gamma$ vanishes, the algorithm terminates.
The following observation is crucial for estimating the time bound.
  \begin{lemma}  \label{lm:break}
Let $q$ be the number of arcs deleted at an iteration with a
breakthrough. Then, excluding the path expansion procedure if applied,
the iteration runs in $O(q)$ time.
  \end{lemma}
  \begin{proof}
Let $\bar W$ be the set of arcs of $\bar\Gamma$ traversed by
TDFS on the iteration, and $W$ the corresponding set in $\Gamma$,
i.e., $W=\{e\in U: e\in \bar W$ or $\bar e\in\bar W\}$.
The iteration runs in $O(q+|W|)$ time, taking into account that each
arc of $P$ not in $Q\cup R$ is contained in $W$. Therefore, it
suffices to show that no arc from $W$ remains in the new graph
$\Gamma$. Suppose this is not so. Then there is a $Z$ to $t$ path $L$
of the old $\Gamma$ that meets $W$ but not $E(Q)\cup E(R)$ (as the
arcs of $Q\cup R$ are deleted).
Let $e=(x,y)$ be the {\em last} arc of $L$ in $W$. Let $b=(y,w)$ be
the next arc of $L$ (it exists since $y=t$ would imply that $e$
is in $A$ but not in $P$, whence $e$ belongs to $Q\cup R$).
Then $b\not\in W\cup E(Q)\cup E(R)$, by the choice of $e$. Two cases
are possible.

(i) $e$ is in $P$. Then $\bar e=(y,x)\in\bar W$. According to
condition (C2), at the time TDFS traversed $\bar e$ from $y$ to $x$
all arcs of $\bar\Gamma$ leaving $y$ had already been traversed.
So $b$ is not in $\bar\Gamma$, implying $b\in E(P)-W$.
Then $b$ is in $Q\cup R$; a contradiction.

(ii) $e$ is not in $P$. Then $e\in \bar W$ and $e$ does not lie on
the final active path $\bar A$ (otherwise $e$ is in $Q\cup R$).
Therefore, TDFS traversed $e$ in both directions. To the time of
traversal of $e$ in backward direction, from $y$ to $x$, all arcs
of $\bar\Gamma$ leaving $y$ have been traversed (at this point the
difference between TDFS and DFS is important). So $b$ is not in
$\bar\Gamma$, whence $b\in E(P)$. Now $\bar b\not\in \bar W$ implies
that $b$ is in $Q\cup R$; a contradiction.
  \end{proof}

\medskip
{\em Non-breakthrough case.} Let $Y$ be the set of nodes visited by
TDFS. Then no arc of $\bar\Gamma$ leaves $Y$. Therefore, in view of
(C1),
  \begin{myitem}
the set of arcs of $\Gamma$ leaving $Y$ consists
of a unique arc $a=(v,w)$, this arc lies on $P$, and the nodes of the
part of $P$ from $z$ to $v$ are contained in $Y$.
  \label{eq:N1}
  \end{myitem}

Since no arc of $\Gamma$ enters $Z$, we also have
  \begin{equation}    \label{eq:N2}
 Y\cap Z=\{z,z'\}.
  \end{equation}

We reduce $\Gamma$ by shrinking its subgraph $\Gamma_Y=(Y,U_Y)$
induced by $Y$ into one node; the formed {\em complex} node $v_Y$ is
identified with $v$.
We call $v$ the {\em root} of $\Gamma_Y$ and store $\Gamma_Y$.
The list of arcs entering $v_Y$ in the new graph is
produced by simply merging the lists $I_x$ for $x\in Y$ from which
the arcs occurring in $\Gamma_Y$ are explicitly removed, using the
lists $O_y$ for $y\in Y$.
(We do not need to correct the outgoing arc lists $O_x$ for
$x\not\in Y$ explicitly, as we explain later.)
Thus, to update $\Gamma$ takes time linear in $|U_Y|$. By~\refeq{N1},
$a$ is the only arc leaving $v_Y$ in $\Gamma$.

From (C1) and~\refeq{N1} it follows that
  \begin{myitem}
the new graph $\Gamma$ is again acyclic, and for each
$x\in Y$, there is a path $P_Y(x)$ from $x$ to the root $v$ in
$\Gamma_Y$.
  \label{eq:N3}
  \end{myitem}

In view of~\refeq{N2}, the set $Z$ is updated as $Z:=Z-\{z,z'\}$.
If there is at least one arc entering $v_Y$, then the new graph
$\Gamma$ and set $Z$ satisfy (C1) and form the input of the next
iteration. The new pre-path $D$ is assigned to be the part of $P$ from
the formed complex node $v_Y$ to $t$. If no arc enters $v_Y$, we
finish the current iteration by removing the nodes and the arcs not
contained in $Z$ to $t$ paths. This further reduce $\Gamma$ and may
reduce $Z$ and shorten $D$.

One can see that the set $U_Y$ is exactly $W$. This and the
construction of pre-paths imply the following.
  \begin{lemma} \label{lm:nobreak}
Let $q$ be the number of arcs deleted by an $i$-th iteration
without a breakthrough.
Then the iteration runs in
$O(q+\max\{0,|D_{i+1}|-|D_i|\})$ time. \qedt
  \end{lemma}

As mentioned above, we do not need to explicitly correct the
outgoing arc lists $O_x$ for $x\not\in Y$
(this would be expensive). Let $\Vscr$ be the current set of all
complex nodes created from the beginning of the algorithm.
We take advantage of the following facts. First, the elements of
$\Vscr$ that are nodes of the current graph (the {\em maximal} complex
nodes) lie on the current pre-path $D$. Second, at an iteration with a
breakthrough, all complex nodes are removed. Third, at an iteration
without a breakthrough, the
subgraph $\Gamma_Y$ forming the new complex node $v_Y$ contains a
subpath of $P$ from its beginning node (by~\refeq{N1}), and the
cleaning procedure (if applied at the iteration) deletes a part of
the updated $P$ from its beginning node as well. Therefore, one can
store $\Vscr$ as a tree in a natural way and use the
{\em disjoint set union} data structure from~\cite{GT-85} to
maintain $\Vscr$. This enables us to efficiently access the head $v_Y$
of any arc $e=(x,v_Y)$ when $e$ is traversed by TDFS (with $O(1)$
amortized time per one arc).

To complete the algorithm description, it remains to explain the
{\em path expansion procedure} to be applied when an iteration with a
breakthrough finds paths $Q,R$ contained complex nodes. It
proceeds in a natural way by recursively expanding complex nodes
occurring in the current $Q,R$ into the corresponding paths $P_Y(x)$ as
in~\refeq{N3} and building $P_Y(x)$ into $Q$ or $R$ (this takes
$O(|P_Y(x)|)$ time). The arc sets of
subgraphs $\Gamma_Y$ extracted during the algorithm are pairwise
disjoint, so the total time for all applications of the procedure
is $O(m)$.

Thus, we can conclude from Lemmas~\ref{lm:break}
and~\ref{lm:nobreak} that the algorithm runs in $O(m)$ time,
yielding Theorem~\ref{tm:bbp}.

\medskip
In the rest of this section we extend the above approach and
algorithm ({\em Algorithm 1}) to a general case of acyclic $(G,u)$.
The auxiliary graph $\Gamma=(X,U)$ and the set $Z$ are constructed
as above, and the capacity $u(e)$ of each arc $e\in U$ is defined
in a natural way. We call an integer $Z$ to
$t$ flow $g$ in $(\Gamma,u)$ {\em balanced} if the flow values out of
``symmetric'' sources are equal, i.e.,
  $$
\div_g(z)=\div_g(\sigma(z))\quad \mbox{for each $z\in Z$,}
  $$
and {\em blocking balanced} if there exists no balanced
flow $g'$ satisfying $g\ne g'\ge g$ (taking into account that $\Gamma$
is acyclic). Then the problem of finding a totally blocking IS-flow in
$(G,u)$ is reduced to {\em problem BBF}: find a balanced blocking flow
for $\Gamma,u,Z,t$.

{\em Algorithm 2} will find a balanced blocking flow
$g$ in the form $g=\alpha_1\chi^{Q_1}+\alpha_1\chi^{R_1}+\ldots+
\alpha_r\chi^{Q_r}+\alpha_r\chi^{R_r}$, where each $\alpha_i$ is a
positive integer, $Q_i$ is a path from some $z\in Z$ to $t$, and $R_i$
is a path from $\sigma(z)$ to $t$.
It iteratively constructs pairs $Q_i,R_i$ for current $\Gamma,u,Z$,
assigns the weight $\alpha_i$ to them as large as possible, and
accordingly reduces the current capacities as $u:=u-
\alpha_i\chi^{Q_i}-\alpha_i\chi^{R_i}$. All arc capacities in $\Gamma$
are positive: once the capacity of an arc becomes
zero, this arc is immediately deleted from $\Gamma$.

Each pair $Q_i,R_i$ is constructed as in Algorithm 1 when
it is applied to the corresponding {\em split-graph} $S=S(\Gamma,u)$.
More precisely (cf.~Section~\ref{sec:theo}),
$S$ is formed by replacing each arc $e=(x,y)$
of $\Gamma$ by two parallel arcs ({\em split-mates}) $e_1,e_2$ from
$x$ to $y$ with the capacities $\lceil u(e)/2\rceil$ and
$\lfloor u(e)/2\rfloor$, respectively. When $u(e)=1$, $e_2$ vanishes in
$S$, and $e_1$ is called {\em critical}. The algorithm maintains $S$
explicitly. The desired pair $Q_i,R_i$ in $(\Gamma,u)$ is
determined by a good pair in $S$ in a natural way.

The main part of an iteration of Algorithm 2 is a slight modification
of an iteration of Algorithm 1. The difference is the following. While
Algorithm 1 deletes {\em all} arcs of the paths $Q,R$ found at an
iteration, Algorithm 2 deletes only a {\em nonempty subset} $B$ of arcs
in $Q\cup R$ (concerning the graph $S$) including all critical arcs in
these paths. One may think that Algorithm 2 essentially treats with
a graph $S$ (ignoring $(\Gamma,u)$) in which
some disjoint pairs of parallel arcs (analogs of split-mates) are
distinguished and the other arcs are regarded as critical, and at each
iteration, the corresponding subset $B\subseteq E(Q)\cup E(R)$ to be
deleted is given by an oracle. Emphasize that the unique arc leaving
a complex node is always critical. Therefore, each complex node in
$Q\cup R$ will be automatically removed. Computing $\alpha_i$'s and
other operations of the algorithm beyond the work with the graph $S$
do not affect the asymptotic time bound.

We now estimate the complexity of an iteration of Algorithm 2. In case
without a breakthrough, properties \refeq{N1},\refeq{N2},\refeq{N3} and
Lemma~\ref{lm:nobreak} (with $S$ instead of $\Gamma$)
remain valid. Note that the arc $a$ in~\refeq{N1} is critical
(since it is a unique arc leaving $Y$); therefore, the arc
leaving the created complex node is critical. Our analysis of the
breakthrough case involves the subset $\bar W_2\subset \bar W$ of
arcs traversed by TDFS in both directions (where $\bar W$ is the
set of all traversed arcs in the corresponding auxiliary graph
$\bar S$). Let $W_2$ be the corresponding set in $S$.
  \begin{lemma} \label{lm:break-gen}
Suppose an iteration of Algorithm 2 results in a breakthrough. Let
$e=(x,y)$ be an arc of $S$ such that $e\in W_2$ or $e\in E(P)\cap W$.
Then any $x$ to $t$ path $L$ in $S$ starting with the arc $e$
contains a critical arc in $Q\cup R$ (and therefore, $e$ vanishes in
the new graph $S$).
  \end{lemma}
 \begin{proof}
Suppose this is not so and consider a counterexample $(e,L)$ with
$|L|$ minimum.
Let $\bar A_0$ be the active path in $\bar S$ just before the
traversal of $e$ or $\bar e$ from $y$ to $x$, and $A_0$ the
corresponding (undirected) path in $S$. At that time, for the set
$\bar O_y$ of arcs of $\bar S$ leaving $y$,
  \begin{myitem}
all arcs in $\bar O_y$ except $\bar e$ (in case $e\in E(P)$) are
already traversed
  \label{eq:ast}
  \end{myitem}
(in view of condition (C2)).
Let $b=(y,w)$ be the second arc of $L$ (existing as
$y=t$ is impossible). By~(\ref{eq:ast}), if $b$ is not in $P$, then
$b\in W$. Also $b\not\in W_2$ and $b\not\in E(P)\cap W$ (otherwise
the part of $L$ from $y$ to $t$ would give a smaller counterexample).
This implies that $b$ belongs to $Q\cup R$, and
therefore, $b$ is not critical. Let $b'$ be the split-mate of $b$.
Considering the path starting with $b'$ and then following
$L$ from $w$ to $t$ (which is smaller than $L$), we similarly
conclude that $b'\not\in W_2$ and $b'\not\in E(P)\cap W$. To come to a
contradiction, we proceed as follows.

The fact that $S$ is acyclic implies that the symmetric difference
(on the arcs) of $P$ and $A_0$ is decomposed into a
path from $Z$ to $t$ and a path from $Z$ to $y$; therefore,
$E(P)\triangle E(A_0)$ contains {\em at most one} arc $a$ leaving $y$.
This and~(\ref{eq:ast}) imply that all arcs in $O_y\cap \bar O_y$
except, possibly, $a$ have been traversed twice; so they are in $W_2$.
Hence, one of $b,b'$ must be in $P$; let for definiteness $b\in E(P)$
(then $b'\not\in E(P)$).

Now $b\not\in W$ implies $b\in E(P)-E(A_0)$, and $b'\in W-W_2$ implies
$b'\in E(A_0)-E(P)$. Thus, both arcs $b,b'$ leaving $y$ are in
$E(P)\triangle E(A_0)$; a contradiction.
  \end{proof}

The running time of an iteration with a breakthrough is
$O(|P|+|W|+q)$, where $q$ is the number of arcs deleted from $S$.
Lemma~\ref{lm:break-gen} allows us to refine this bound as
$O(|Q|+|R|+q)$. Combining this with Lemma~\ref{lm:nobreak}, we can
conclude that, up to a constant factor, the total time of Algorithm 2
is bounded from above by $m$ plus the sum $\Sigma$ of lengths of paths
$Q_1,R_1,\ldots,Q_r,R_r$ in the representation of the flow $g$
constructed by the algorithm. Since $|Q_i|,|R_i|\le n$ and $r\le 2m$
(as each iteration decreases the arc set of $S$), $\Sigma$ is $O(nm)$.
Also $\Sigma$ does not exceed the sum of the transit capacities $u(x)$
of inner nodes $x$ of $\Gamma$ (assuming, without loss of generality,
that no arc goes from $s$ to $s'$). Thus, Theorem~\ref{tm:ac-unit} is
generalized as follows.
  \begin{theorem}  \label{tm:ac-gen}
For an acyclic capacitated skew-symmetric network $N$ with $O(n)$ nodes
and $O(m)$ arcs, a totally blocking IS-flow can be found in
$O(\min\{m+\Delta(N),nm\})$ time.
  \end{theorem}

Together with Corollary~\ref{cor:cvb} and Lemma~\ref{lm:root},
this yields the desired generalization.
\begin{corollary} \label{cor:gen_cap}
SBFM can be implemented so that it finds a maximum IS-flow in an
arbitrary skew-symmetric network $N$ in
$O(\min\{n^2m,\sqrt{\Delta(N)}(m+\Delta(N))\})$ time.
  \end{corollary}

%-------  SECTION 10 --------------

\section{\Large Applications to Matchings}\label{sec:mat}

Apply the reduction of the maximum u-capacitated b-matching problem
(CBMP) in a graph $G'=(V',E')$ to the maximum IS-flow problem in
a network $N=(G=(V,E),\sigma,u,s)$; see Section~\ref{sec:back}.
The best time bound for a general case of CBMP is attained by applying
the algorithm of Section~\ref{sec:gisa}. Theorem~\ref{tm:gisa} implies
the following.
  \begin{corollary} \label{cor:bmat}
CBMP can be solved in $O(M(n,m)+nm)$ time, where $n:=|V'|$ and
$m:=|E'|$.
   \end{corollary}
When the input functions $u,b$ in CBMP are small enough, the
transit capacities of nodes in $N$ become small as well. Then the
application of the shortest blocking IS-flow method may result in a
competitive or even faster algorithm for CBMP.
Let the capacities of all edges of $G'$ be ones. We have
$\Delta(N)=O(m)$ in general, and $\Delta(N)=O(n)$ if $b$ is all-unit.
Then Corollary~\ref{cor:unit} yields the same time bounds as
in~\cite{gab-83,MV-80} for the corresponding cases.
\begin{corollary} \label{cor:bmat_un}
SBFM (with the fast implementation of a phase as in
Section~\ref{sec:acyc}) solves the maximum degree-constrained subgraph
(or b-factor) problem in $O(m^{3/2})$ time and solves the maximum
matching problem in a general graph in $O(\sqrt{n}m)$ time.
  \end{corollary}

Feder and Motwani~\cite{FM-91} elaborated a clique compression
technique and used it to improve the $O(\sqrt{n}m)$ bound for the
maximum bipartite matching problem to $O(\sqrt{n}m\log(n^2/m)/\log{n})$.
We explain how to apply a similar approach to a special case of MSFP,
lowering the bound for dense nonbipartite graphs. We need a brief
review of the method in~\cite{FM-91}.

Let $H=(X,Y,E)$ be a bipartite digraph, where $E\subseteq X\times Y$,
$|X|=|Y|=n$ and $|E|=m$. A {\em (bipartite) clique} is a complete
bipartite subgraph $(A,B,A\times B)$ of $H$, denoted by $C(A,B)$.
Define the {\em size} $s(C)$ of $C=C(A,B)$ to be $|A|+|B|$.
A {\em clique partition} of $H$ is a collection
$\Cscr$ of cliques whose arc sets form a partition of $E$; the {\em
size} $s(\Cscr)$ of $\Cscr$ is the sum of sizes of its members.

Let a constant $0<\delta<1/2$ be fixed. Then a clique $C(A,B)$ of $H$
is called a $\delta$-{\em clique} if $|A|=\lceil n^{1-\delta}\rceil$
and $|B|=\lfloor \delta\log{n}/\log(2n^2/m) \rfloor$. It is shown
in~\cite{FM-91} that a $\delta$-clique exists.

The {\em clique partition algorithm} in~\cite{FM-91}
finds a $\delta$-clique $C_1$ in the initial graph $H_1=(X,Y,E=:E_1)$
and deletes the arcs of $C_1$, obtaining the next graph $H_2=
(X,Y,E_2)$. Then it finds a $\delta$-clique $C_2$
(concerning the number of arcs of $H_2$) and delete the arcs of $C_2$
from $H_2$, and so on while the number of arcs of the current graph
is at least $2n^{2-\delta}$ and the $Y$-part of a $\delta$-clique is
nonempty. The remaining arcs are partitioned into
cliques consisting of a single arc each. So the cliques $C_i$ extracted
during the algorithm form a clique partition. The running time of the
algorithm is estimated as the sum of bounds $\tau(C_i)$ on the time to
extract the cliques $C_i$ plus a time bound $\tau'$ to maintain a
certain data structure (so-called neighborhood trees). One shows that
\begin{myitem}
the algorithm runs in $O(\sqrt{n}m\beta)$ time and finds a clique
partition $\Cscr$ of $H$ such that $s(\Cscr)=O(m\beta)$, where
$\beta:=\log(n^2/m)/\log{n}$.
  \label{eq:part}
  \end{myitem}

Suppose we wish to find a maximum matching in a bipartite graph or,
equivalently, to find a maximum integer flow from $s$ to $t$ in a
digraph $G$ with unit arc capacities, node set $X\cup Y\cup\{s,t\}$
and arc set $E\cup(s\times X)\cup(Y\times t)$, where $E\subseteq
X\times Y$. One may assume $|X|=|Y|=n$. Using the above algorithm,
form a clique partition $\Cscr$ as in~\refeq{part} for $(X,Y,E)$.
Transform each clique $C(A,B)$ in $\Cscr$ into a star by replacing
its arcs by a node $z$, arcs $(x,z)$ for all $x\in A$ and arcs
$(z,y)$ for all $y\in B$.
There is a natural one-to-one correspondence between the
$s$ to $t$ paths in $G$ and those in the resulting graph $G^\ast$,
and the problem for $G^\ast$ is equivalent to that for $G$.
Compared with $G$, the graph $G^\ast$ has $|\Cscr|$ additional nodes
but the number $m^\ast$ of its arcs becomes
$2n+s(\Cscr)$, or $O(m\beta)$. Given a flow in $G^\ast$, any (simple)
augmenting path of length $q$ meets exactly $(q-1)/2$ nodes in
$X\cup Y$, and these nodes have unit transit capacities. This implies
that Dinits' algorithm has $O(\sqrt{n})$ phases (arguing as
in~\cite{ET-75,kar-73-2}). Since each phase takes $O(m^\ast)$ time, the
whole algorithm runs in $O(\sqrt{n}m\beta)$ time, as desired.

Now suppose $H=(X,Y,E,\sigma)$ is a skew-symmetric bipartite graph
without parallel arcs, where the sets $X$ and $Y$ are symmetric
each other. We modify the above method as follows. Note that any two
symmetric cliques in $H$ are disjoint (otherwise some
$x\in X$ is adjacent to $\sigma(x)$, implying the existence of two arcs
from $x$ to $\sigma(x)$). We call a clique partition $\Cscr$
{\em symmetric} if $C\in\Cscr$ implies $\sigma(C)\in\Cscr$. An
iteration of the {\em symmetric clique partition algorithm} works as in
the previous algorithm, searching for a $\delta$-clique $C'$ in the
current $H$, but then deletes the arcs of {\em both} $C'$ and
$\sigma(C')$. Let the algorithm construct a partition $\Cscr'$
consisting of cliques $C'_1,\sigma(C'_1),\ldots,C'_r,\sigma(C'_r)$
obtained in this order.

To estimate the size of $\Cscr'$ and the running time, imagine we
would apply the previous algorithm to our $H$ (ignoring the
fact that $H$ is skew-symmetric). Let the resulting partition $\Cscr$
be formed by cliques $C_1,\ldots,C_q$ (in this order).
Note that for a bipartite graph with $n$ nodes and $m$ arcs, both the
number $e(C)$ of arcs of a $\delta$-clique $C$ and its size $s(C)$
are computed uniquely, and these are monotone functions in $m$,
as well as the above-mentioned time bound $\tau(C)$ (indicated
in~\cite{FM-91}). Moreover, one can check that $m'\le m$ and
$e(C')\ne 0$ imply $m'-2e(C')\le m-e(C)$,
where $C'$ is a $\delta$-clique in a graph with
$n$ nodes and $m'$ arcs. Using these, we can conclude that
$r\le q$ and that for $i=1,\ldots r$, $s(C'_i)\le s(C_i)$ and
$\tau(C'_i)\le \tau(C_i)$. Then $s(\Cscr')\le 2s(\Cscr)$,
implying $s(\Cscr)=O(m\beta)$, by~\refeq{part}.
Also the time of the modified algorithm is $O(\sqrt{n}m\beta)$
(by~\refeq{part} and by the fact that the above
bound $\tau'$ remains the same) plus the time needed to treat the
symmetric cliques $\sigma(C'_i)$, which is $O(m)$.

Finally, the graph $H^\ast$ obtained from $H$ by transforming the
cliques $C'_1,\sigma(C'_1),\ldots, C'_r,\sigma(C'_r)$ into stars
has a naturally induced skew-symmetry. By the above argument, $H^\ast$
has $O(m\beta)$ arcs, and computing $H^\ast$ takes
$O(\sqrt{n}m\beta)$ time. Apply such a transformation to the input
graph of MSFP arising from an instance of the maximum matching problem.
Arguing as in the bipartite matching case above and as in the proof of
Lemma~\ref{lm:root}, we conclude with the following.
  \begin{theorem} \label{tm:compr}
A maximum matching in a general graph with $n$ nodes and $m$
edges can be found in $O(\sqrt{n}m\log(n^2/m)/\log{n})$ time.
  \end{theorem}

{\bf Acknowledgements.}
We thank the anonymous referees for suggesting improvements in the
original version of this paper and pointing out to us important
references.


\small
\begin{thebibliography}{99}

\bibitem{ans-85}
R.P.~Anstee. An algorithmic proof of Tutte's $f$-factor theorem.
{\sl J. of Algorithms} {\bf 6} (1985) 112--131.

\bibitem{ans-87}
R.P.~Anstee. A polynomial algorithm for $b$-matchings: An alternative
approach. {\sl Information Proc. Letters} {\bf 24} (1987) 153--157.

\bibitem{blu-90}
N.~Blum.
A new approach to maximum matching in general graphs.
In {\sl Automata, Languages and Rpogramming} [Lecture Notes in Comput.
Sci. 443] (Springer, Berlin, 1990), pp. 586--597.

\bibitem{din-70}
E.~A. Dinic.
Algorithm for solution of a problem of maximum flow in networks with
power estimation.
{\sl Soviet Math. Dokl.} {\bf 11} (1970) 1277--1280.

\bibitem{edm-65}
J.~Edmonds.
Paths, trees and flowers.
{\sl Canada J. Math.} {\bf 17} (1965) 449--467.

\bibitem{EJ-70}
J.~Edmonds and E.~L. Johnson.
Matching, a well-solved class of integer linear programs.
In R.~Guy, H.~Haneni, and J.~Sch{\"o}nhein, eds,
{\sl Combinatorial Structures and Their Applications},
Gordon and Breach, NY, 1970, pp. 89--92.

\bibitem{EK-72}
J.~Edmonds and R.~M. Karp.
Theoretical improvements in algorithmic efficiency for network flow
problems.
{\sl J. Assoc. Comput. Mach.} {\bf 19} (1972) 248--264.

\bibitem{ET-75}
S.~Even and R.~E. Tarjan.
Network flow and testing graph connectivity.
{\sl SIAM J. Comput.} {\bf 4} (1975) 507--518.

\bibitem{FM-91}
T.~Feder and R.~Motwani.
Clique partitions, graph compression and speeding-up algorithms.
In {\sl Proc. 23rd Annual ACM Symp. on Theory of Computing},
1991, pp.~123--133.

\bibitem{FF-62}
L.~R. Ford and D.~R. Fulkerson.
{\sl Flows in Networks}.
Princeton Univ. Press, Princeton, NJ, 1962.

\bibitem{FJ-99}
C.~Fremuth-Paeger and D.~Jungnickel.
Balanced network flows. Parts I--III.
{\sl Networks} {\bf 33} (1999) 1--56.

\bibitem{gab-83}
H.~N. Gabow.
An efficient reduction technique for degree-constrained subgraph and
bidirected network flow problems.
{\sl Proc. of STOC} {\bf 15} (1983) 448--456.

\bibitem{GT-85}
H.~N. Gabow and R.~E. Tarjan.
A linear-time algorithm for a special case of disjoint set union.
{\sl J. Comp. and Syst. Sci.} {\bf 30} (1985) 209--221.

\bibitem{GT-91}
H.~N. Gabow and R.~E. Tarjan.
Faster scaling algorithms for general graph-matching problems.
{\sl J. ACM} {\bf 38} (1991) 815--853.

\bibitem{GK-95}
A.~V. Goldberg and A.~V. Karzanov.
Maximum skew-symmetric flows.
In P.~Spirakis, ed., {\sl Algorithms -- ESA '95} (Proc. 3rd European
Symp. on Algorithms),  {\sl Lecture Notes in Computer Sci.}
{\bf 979}, 1995, pp. 155--170.

\bibitem{GK-96}
A.~V. Goldberg and A.~V. Karzanov.
Path problems in skew-symmetric graphs.
{\sl Combinatorica} {\bf 16} (1996) 129--174.

\bibitem{GK-99}
A.~V. Goldberg and A.~V. Karzanov.
Maximum skew-symmetric flows and their applications to b-matchings.
{\sl Preprint} 99-043, SFB 343, Bielefeld Universit\"at, Bielefeld,
1999, 25 pp.

\bibitem{HK-73}
J.~E. Hopcroft and R.~M. Karp.
An $n^{5/2}$ algorithm for maximum matching in bipartite graphs.
{\sl SIAM J. Comput.} {\bf 2} (1973) 225--231.

\bibitem{kar-70}
A.~V. Karzanov.
\`Ekonomny\u{i} algoritm nakhozhdeniya bikomponent grafa
[Russian; An efficient algorithm for finding the bicomponents of a
graph].
In {\sl Trudy Tret'e\u{i} Zimne\u{i} Shkoly po Matematicheskomu
Programmirovaniyu i Smezhnym Voprosam} [Proc. of 3rd Winter School on
Mathematical Programming and Related Topics], issue {\bf 2}. Moscow
Engineering and Construction Inst. Press, Moscow, 1970, pp.~343-347.

\bibitem{kar-73}
A.V.~Karzanov.
Tochnaya otsenka algoritma nakhozhdeniya maksimal'nogo potoka,
primenennogo k zadache ``o predstavitelyakh''
[Russian;  An exact estimate of an algorithm for
finding a maximum flow, applied to the problem ``of representatives''].
In {\sl Voprosy Kibernetiki} [Problems of Cybernetics],
volume~{\bf 3}. Sovetskoe Radio, Moscow, 1973, pp.~66--70.

\bibitem{kar-73-2}
A.V.~Karzanov.
O nakhozhdenii maksimal'nogo potoka v setyakh spetsial'nogo vida i
nekotorykh prilozheniyakh
[Russian;  On finding maximum flows in networks
with special structure and some applications].
In {\sl Matematicheskie Voprosy Upravleniya Proizvodstvom}
[Mathematical Problems for Production Control],
volume~{\bf 5}.
Moscow State University Press, Moscow, 1973, pp. 81-94.

\bibitem{KS-93}
W.~Kocay and D.~Stone. Balanced network flows.
{\sl Bulletin of the ICA} {\bf 7} (1993) 17--32.

\bibitem{law-76}
E.~L. Lawler.
{\sl Combinatorial Optimization: Networks and Matroids}.
Holt, Reinhart, and Winston, New York, NY., 1976.

\bibitem{LP-86}
L.~Lov{\'a}sz and M.~D. Plummer.
{\sl {Matching Theory}}.
Akad{\'e}miai Kiad{\'o}, Budapest, 1986.

\bibitem{MV-80}
S.~Micali and V.V.~Vazirani.
An $O(\sqrt{V} E)$ algorithm for finding maximum matching in general
graphs.
{\sl Proc. of the 21st Annual IEEE Symposium in Foundation of Computer
Science}, 1980, pp.~71--109.

\bibitem{sch-03}
A.~Schrijver. {\sl Combinatorial Optimization. Polyhedra and
Efficiency}, Volume~A. ({\sl Algorithms and Combinatorics} {\bf 24}),
Springer, Berlin and etc., 2003.

\bibitem{tar-72}
 R.~Tarjan. Depth-first search and linear graph algorithms.
{\sl SIAM J. Comput.} {\bf 1} (1972) 146--160.

\bibitem{tut-67}
W.T.~Tutte.
Antisymmetrical digraphs.
{\sl Canadian J. Math.} {\bf 19} (1967) 1101--1117.

\bibitem{vaz-90}
V.V.~Vazirani.
A theory of alternating paths and blossoms for proving correctness of
the $O(\sqrt{V} E)$ general graph maximum matching algorithm.
In: R.~Kannan and W.R.~Cunningham, eds., {\sl Integer Programming and
Combinatorial Optimization} (Proc. 1st IPCO Conference), University of
Waterloo Press, Waterloo, Ontario, 1990, pp.~509--535.

\end{thebibliography}
\end{document}